\documentclass[12pt]{article}
\usepackage{graphicx,psfrag,epsfig}
\usepackage{amssymb,amsmath,amscd,amsthm}
\usepackage{graphicx,psfrag,epsfig}

\usepackage{graphicx}
\usepackage[active]{srcltx}

\newtheorem{theorem}{Theorem}[section]
\newtheorem{corollary}[theorem]{Corollary}

\newtheorem{lemma}[theorem]{Lemma}

\date{}

\begin{document}

\date{}
\title{Nonlinear Stochastic Perturbations of Dynamical Systems and
 Quasi-linear Parabolic PDE's with a Small Parameter}
\author{ M. Freidlin\footnote{Dept of Mathematics, University of Maryland,
College Park, MD 20742, mif@math.umd.edu}, L.
Koralov\footnote{Dept of Mathematics, University of Maryland,
College Park, MD 20742, koralov@math.umd.edu}
} \maketitle

\begin{abstract}
In this paper we describe the asymptotic behavior, in the
exponential time scale, of  solutions to quasi-linear parabolic
equations with a small parameter at the second order term and the
long time behavior of corresponding diffusion processes. In
particular, we discuss the exit problem and metastability for the
processes corresponding to quasi-linear initial-boundary value
problems.
\end{abstract}



\section{Introduction}

Consider a dynamical system
\begin{equation} \label{dsyst}
\dot{X}^{x}_t = b( X^{x}_t),~~X^x_0 = x \in \mathbb{R}^d,
\end{equation}
together with its stochastic perturbations
\begin{equation} \label{perturb}
d X^{x,\varepsilon}_t = b(  X^{x,\varepsilon}_t) d t+\varepsilon
\sigma ( X^{x,\varepsilon}_t) d W_t,~~X^{x,\varepsilon}_0 = x \in
\mathbb{R}^d.
\end{equation}
Here $\varepsilon > 0$ is a small parameter, $W_t$ is a Wiener
process in $ \mathbb{R}^d$, and the coefficients $\sigma$ and~$b$
are assumed to be Lipschitz continuous. The diffusion matrix $a(x)
= (a_{ij}(x)) = \sigma(x) \sigma^*(x)$ is assumed to be
non-degenerate for all $x$.

Let $D$ be a bounded domain in $ \mathbb{R}^d$ with infinitely
smooth boundary $\partial D$. In this paper, with the exception of
the last section, we assume that there is a point $x_0 \in D$ such
that for each $x \in D$ the trajectory of the dynamical system
(\ref{dsyst}) starting at $x$ is attracted to $x_0$. We assume
that ${(b(x), n(x)) < 0}$ for $x \in \partial D$, where $n(x)$ is
the exterior normal to the boundary of $D$. Let $\tau^\varepsilon
= \min\{t: X^{x,\varepsilon}_t \in \partial D\}$ be the first time
when  $ X^{x,\varepsilon}_t$ reaches the boundary of $D$.

If $\varepsilon$ is small, then on any finite time interval the
trajectories of the process $ X^{x,\varepsilon}_t$ defined by
(\ref{perturb}) are close to the corresponding non-perturbed
trajectory with probability close to one. Therefore, with high
probability $X^{x,\varepsilon}_t$ enters a small neighborhood of
the equilibrium point $x_0$ before leaving $D$. The process
eventually exits $D$ as a result of large deviations of $
X^{x,\varepsilon}_t$ from $X^x_t$ (\cite{FW}, see also \cite{OV}).
The large deviations are governed by the normalized action
functional
\[
S_{0, T}(\varphi) = \frac{1}{2}\int_{0}^{T} \sum_{i,j = 1}^d
a^{ij}(\varphi_t)(\dot{\varphi}^i_t -
b_i(\varphi_t))(\dot{\varphi}^j_t - b_j(\varphi_t)) d t,~~T \geq
0,~ \varphi \in C([0,T], \overline{D}),
\]
and the quasi-potential
\[
V(x_0,x) =  \inf_{T, \varphi} \{ S_{0,T}(\varphi): \varphi \in
C([0,T], \overline{D}), \varphi(0) = x_0, \varphi(T) = x \},~~x
\in \overline{D}.
\]
Here $a^{ij}$ be the elements of the inverse matrix, that is
$a^{ij} = (a^{-1})_{ij}$, and $S_{0, T}(\varphi) = +\infty$ if
$\varphi$ is not absolutely continuous. It is proved in
$\cite{FW}$ that $\varepsilon^2 \ln \tau^\varepsilon$ converges in
probability, as $\varepsilon \downarrow 0$, to $V_0 = \min_{x \in
\partial D} V(x_0,x)$. Moreover, if the minimum $V_0$ of $V(x_0,x)$ on
$\partial D$ is achieved at a unique point $x^*$ (which is true in
the generic case), then $ X^{x,\varepsilon}_{\tau^\varepsilon}$
converges to $x^*$ in probability as~${\varepsilon \downarrow 0}$.

These statements imply various results for PDE's with a small
parameter at the second order derivatives. In particular, consider
the following initial-boundary value problem:
\begin{equation} \label{y1}
\frac{\partial w^\varepsilon(t,x)}{\partial t}  =
\frac{\varepsilon^2}{2} \sum_{i,j=1}^d a_{ij}(x) \frac{ \partial^2
w^\varepsilon (t,x)}{\partial x_i \partial x_j} + b(x) \cdot
\nabla_x w^\varepsilon (t,x),~~x \in D,~~t
> 0,
\end{equation}
\begin{equation} \label{y2}
w^\varepsilon(0,x) = g(x),~x \in D,~~~~~w^\varepsilon(t,x) =
g(x),~~t \geq 0,~x
\in
\partial D,
\end{equation}
where $g$, for the sake of brevity, is assumed to be continuous on
$\overline{D}$. The case when $w^\varepsilon(t,x)|_{x \in \partial
D} = \psi(x)$ with $\psi \neq g$  can be considered in a similar
way. Assume that the minimum $V_0$ of $V(x_0,x)$ on $\partial D$
is achieved at a unique point $x^*$. Let $t: \mathbb{R}^+
\rightarrow \mathbb{R}$ be a function such that $t(\varepsilon)
\asymp \exp(\lambda/\varepsilon^2)$ as $\varepsilon \downarrow 0$
with $\lambda > 0$, that is $\ln (t(\varepsilon)) \sim
\lambda/\varepsilon^2$ as $\varepsilon \downarrow 0$. Then
\[
 \lim_{\varepsilon \downarrow
0} w^\varepsilon(t(\varepsilon), x) = g(x_0),~~{\rm if}~\lambda <
V_0,
\]
\[
\lim_{\varepsilon \downarrow 0} w^\varepsilon(t(\varepsilon), x) =
g(x^*),~~{\rm if}~\lambda > V_0,
\]
for $x \in D$.

Note that the solution to (\ref{y1})-(\ref{y2}) can be expressed
in terms of the transition semigroup associated with the family of
processes $X^{x,\varepsilon}_t$, $x \in \overline{D}$. Namely, let
$T^\varepsilon_t g(x) = \mathrm{E}g(X^{x,\varepsilon}_{t \wedge
\tau^\varepsilon})$, $g \in C(\overline{ D})$. Then the function
$w^\varepsilon(t,x) = T^\varepsilon_t g(x)$ is the solution to
(\ref{y1})-(\ref{y2}). The semigroup $T^\varepsilon_t$ can be
viewed as a small perturbation of the semigroup of shifts $T_t
g(x) = g(X^{x}_{t})$ associated with the dynamical system
(\ref{dsyst}).

More general perturbations of $T_t$ may lead to nonlinear
semigroups. Namely, consider the following problem:
\begin{equation} \label{eq1}
\frac{\partial u^\varepsilon(t,x)}{\partial t}  = L^\varepsilon
u^\varepsilon := \frac{\varepsilon^2}{2} \sum_{i,j=1}^d
a_{ij}(x,u^\varepsilon) \frac{ \partial^2 u^\varepsilon
(t,x)}{\partial x_i \partial x_j} + b(x) \cdot \nabla_x
u^\varepsilon (t,x),~~x \in D,~~t
> 0,
\end{equation}
\begin{equation} \label{bc}
u^\varepsilon(0,x) = g(x),~x \in D;~~~~~u^\varepsilon(t,x) =
g(x),~~t \geq 0,~x
\in
\partial D.
\end{equation}
When the coefficients are sufficiently smooth and the matrix $a$
is positive-definite, the solution $u^\varepsilon$ exists and is
unique in the appropriate function space (see Section~\ref{qle}
below). We can now define the semigroup $T^\varepsilon_t$ on
$C(\overline{D})$ via $T^\varepsilon_tg(x) = u^\varepsilon(t,x)$,
where $u^\varepsilon$ is the solution of (\ref{eq1})-(\ref{bc})
with initial-boundary data $g$.

For $t > 0$ and $x \in \overline{D}$, we can define
$X^{t,x,\varepsilon}_s$, $s \in [0,t]$, as the process which
starts at $x$ and solves
\begin{equation} \label{smttt}
d X^{t,x,\varepsilon}_s = b( X^{t,x,\varepsilon}_s) d s +
\varepsilon \sigma(X^{t,x,\varepsilon}_s, u^\varepsilon(t-s,
X^{t,x,\varepsilon}_s)) d W_{s},~~s \leq \tau^\varepsilon  \wedge
t,
\end{equation}
\[
\tau^\varepsilon = \min\{s: X^{t,x,\varepsilon}_s \in
\partial D\},~~~~~~~~ X^{t,x,\varepsilon}_s =
 X^{t,x,\varepsilon}_{\tau^\varepsilon},~~\tau^\varepsilon \leq s
 \leq t,
\]
where $\sigma_{ij}$, $1 \leq i,j \leq d$, are Lipschitz continuous
and such that $\sigma \sigma^* = a$.
The process $X^{t,x,\varepsilon}_s$ will be called the nonlinear
stochastic perturbation of (\ref{dsyst}). More precisely,
$X^{t,x,\varepsilon}_s$ corresponds to the nonlinear semigroup
defined by (\ref{eq1}). As in the linear case, we have the
following relation between $u^\varepsilon$ and the process
$X^{t,x,\varepsilon}_s$:
\[
u^\varepsilon(t,x) = \mathrm{E}g(X^{t,x,\varepsilon}_{t \wedge
\tau^\varepsilon}).
\]

One of the important questions in the study of parabolic linear
and nonlinear equations is the one concerning the behavior of
solutions (or, in probabilistic terms, behavior of the
corresponding diffusion process) as $t \rightarrow \infty$. In our
case, when the small parameter $\varepsilon^2$ is present in front
of the second order term, the limit of $u^\varepsilon(t,x)$ as
$\varepsilon \rightarrow 0$, $t \rightarrow \infty$, depends on
the manner in which $(\varepsilon, t)$ approaches $(0, \infty)$.
In the linear case this problem has been studied in \cite{F77}
(see also \cite{FW}, \cite{OV}).

In Section~\ref{asc} we study the asymptotic behavior of solutions
to (\ref{eq1})-(\ref{bc}) when $\varepsilon \downarrow 0$ and $t =
t(\varepsilon) \asymp \exp(\lambda/\varepsilon^2)$. As a first
step, we shall introduce a family of linear problems which can be
obtained from (\ref{eq1})-(\ref{bc}) by replacing the second
variable in the coefficients $a_{ij}$ in the right hand side of
(\ref{eq1}) by a constant $c$. The asymptotics of $u^\varepsilon$
can be then expressed in terms of the functions $V_0(c)$ and
$g(x^*(c))$, where $V_0(c)$ is the minimum of the quasi-potential
of the linear problem and $x^*(c)$ is the point where this minimum
is achieved.

In Section~\ref{exitprob} we study the exit problem for the
process $X^{t,x,\varepsilon}_s$. We shall see that new effects
appear when nonlinear stochastic perturbations are considered. In
particular, even in the generic case, the distribution of the exit
location $X^{t,x,\varepsilon}_{t \wedge \tau^\varepsilon}$ need
not be concentrated in one point.


Some related problems concern the notion of metastability for
nonlinear perturbations of dynamical systems with several
equilibrium points. Let us consider the dynamical system
(\ref{dsyst}) in $ \mathbb{R}^d$ and its perturbations
(\ref{perturb}). As before, $\sigma$ and~$b$ are assumed to be
Lipschitz continuous. Now we shall assume that the system has a
finite number of asymptotically stable equilibrium points
$x_1,...,x_k$ such that for almost every $x \in \mathbb{R}^d$,
with respect to the Lebesgue measure, the trajectory  of
(\ref{dsyst}) starting at $x$ is attracted to one of the points
$x_1,...x_k$.
 We shall also assume that the vector
field $b$ satisfies $(b(x), x) \leq A - B|x|^2$ for some positive
constants $A$ and $B$. The case of more general asymptotically
stable attractors (for instance, limit cycles) can be considered
similarly, however for the sake of brevity we restrict ourselves
to the case of equilibriums.

The general theory of metastability was developed in \cite{F77} in
the framework of large deviations (see also \cite{FW},
\cite{FPhD}, \cite{OV}). It was shown, in particular, that for a
generic vector field $b$ satisfying the assumptions above, for
almost every $x \in \mathbb{R}^d$ and $\lambda > 0$, with
probability which tends to one when $\varepsilon \downarrow 0$,
the trajectory $X^{x,\varepsilon}_t$ of (\ref{perturb}) spends
most of the time in the time interval $[0,
\exp(\lambda/\varepsilon^2)]$ near a point $x^\lambda \in
\{x_1,...x_k\}$. This point is called the metastable state for the
trajectory starting at~$x$ in the time scale
$\exp(\lambda/\varepsilon^2)$. The metastable state can be
determined by examining the values of the quasi-potential. Namely,
let
\[
V_{ij} = V(x_i,x_j) =  \inf_{T, \varphi} \{ S_{0,T}(\varphi):
\varphi \in C([0,T], \mathbb{R}^d), \varphi(0) = x_i, \varphi(T) =
x_j \},~~1 \leq i,j \leq m.
\]
These numbers determine a hierarchy of cycles along which the
system switches from one metastable state to another with the
growth of $\lambda$ (\cite{F77}).

We can also study metastability for nonlinear perturbations of
dynamical systems. It turns out that now the transition between
the equilibrium points does not occur ``immediately in the
exponential time scale''. This implies that now metastable states
should be replaced by metastable distributions between the
equilibriums. The description of metastable distributions is based
on the study of the asymptotic behavior of solutions to
(\ref{eq1})-(\ref{bc}) when $\varepsilon \downarrow 0$ and $t =
t(\varepsilon) \asymp \exp(\lambda/\varepsilon^2)$. Note that
metastable distributions also arise in \cite{AF}, \cite{F}, but
for reasons which are different from what is discussed in this
paper. Such a modification to the notion of metastability leads to
a modified notion of stochastic resonance.

We briefly address the problems of metastability in
Section~\ref{genex}, where we also consider other generalizations
and some examples. The issue of metastability in the case of an
arbitrary number of equilibrium points and cycles will be
addressed in a forthcoming paper.

In this paper we considered nonlinear perturbations of a system
with an asymptotically stable equilibrium. In this case, the exit
from a domain containing this equilibrium occurs due to large
deviations, and the exit time and exit distribution essentially
depend on the perturbation. A related singular perturbation
problem arises in the case when the equilbrium is stable but not
asymptotically stable: for instance when the unperturbed system is
Hamiltonian. Nonlinear stochastic perturbations in this case lead
to a nonlinear version of the averaging principle. Say, in the
case of one degree of freedom, the limiting slow motion is a
diffusion process corresponding to a nonlinear operator on the
graph (compare with \cite{FW}, Chapter 8) related to the
Hamiltonian. We will consider these problems in one of the
forthcoming papers.

\section{Preliminaries and Notations}
\subsection{Quasi-Linear Equation} \label{qle}

Let $D \subset \mathbb{R}^d$ be a bounded domain with infinitely
smooth boundary $\partial D$. We shall say that $f: D \rightarrow
\mathbb{R}$ belongs to $C^2(D)$ if $f$ and all of its partial
derivatives up to the second order are bounded and continuous in
$D$.  We shall say that a function $f: (0,\infty) \times D
\rightarrow \mathbb{R}$ belongs to $C^{1,2}((0,\infty) \times D)$
if $f$, its partial derivative in $t$, and all of its partial
derivatives up to the second order in $x$ are bounded and
continuous in $(0,\infty) \times D$. Note that a function $f \in
C^2(D)$ can be extended to a continuous function on $\overline{D}$
and $f \in C^{1,2}((0,\infty) \times D)$  can be extended to a
continuous function on $ [0,\infty) \times \overline{D}$.

Let $a_{ij} = a_{ji} \in C^2 (D \times \mathbb{R})$, $1 \leq i,j
\leq d$, and $b_i \in C^2(D)$, $1 \leq i \leq d$. We also assume
that there is a positive constant $k$ such that $k|\xi|^2 \leq
\sum_{i,j =1}^d a_{ij}(x, u) \xi_i \xi_j$, $x \in D, u \in
\mathbb{R}$, $\xi \in \mathbb{R}^d$. Let $g$ be an infinitely
smooth function defined in a neighborhood of $\overline{D}$.

If  $\varepsilon > 0$ and the coefficients $a$ and $b$ and the
function $g$ satisfy the assumptions listed above, then the
equation (\ref{eq1})-(\ref{bc}) has a unique solution in the class
of functions $C^{1,2}((0,\infty) \times D) \cap C( [0,\infty)
\times \overline{D})$ (see Theorem 5, Chapter~6.2 of \cite{K}). If
$g$ were to be only continuous on $\overline{D}$, the existence
and uniqueness of solutions to (\ref{eq1})-(\ref{bc}) would hold
in the class of functions which are locally $C^{1,2}$-smooth
inside $(0,\infty) \times D$ and continuous up to the boundary.
However, to simplify notations in later sections, we impose the
smoothness condition on $g$.

\subsection{Action Functional} \label{af}


 Let $\alpha$ be a symmetric $d \times d$ matrix whose elements
$\alpha_{ij}$ are bounded and Lipschitz continuous on $
\mathbb{R}^d$ and  satisfy $k|\xi|^2 \leq \sum_{i,j =1}^d
\alpha_{ij}(x) \xi_i \xi_j$, $x \in \mathbb{R}^d$, $\xi \in
\mathbb{R}^d$.  Let $\alpha^{ij}$ be the elements of the inverse
matrix, that is $\alpha^{ij} = (\alpha^{-1})_{ij}$, and $\sigma$
be a square matrix such that $\alpha = \sigma \sigma^*$. We choose
$\sigma$ in such a way that $\sigma_{ij}$ are also bounded and
Lipschitz continuous.

Let $S^\alpha_{0, T}$ be the normalized action functional for the
family of processes $X^{x,\varepsilon}_t$ satisfying
\[
d X^{x,\varepsilon}_t = b(  X^{x,\varepsilon}_t) d t + \varepsilon
\sigma ( X^{x,\varepsilon}_t) d W_t,
\]
 where $b$ is  a bounded Lipschitz continuous vector field
on $ \mathbb{R}^d$. Thus
\[
S^\alpha_{0, T}(\varphi) = \frac{1}{2}\int_{0}^{T} \sum_{i,j =
1}^d \alpha^{ij}(\varphi_t)(\dot{\varphi}^i_t -
b_i(\varphi_t))(\dot{\varphi}^j_t - b_j(\varphi_t)) d t
\]
for absolutely continuous $\varphi$ defined on $[0, T]$,
$\varphi_0 = x$,  and $ S^\alpha_{0, T}(\varphi) = \infty$ if
$\varphi$ is not absolutely continuous or if $\varphi_0 \neq x$
(see \cite{FW}). Let $V^\alpha(x,y)$ be the quasi-potential for
the family $X^{x,\varepsilon}_t$ in $\overline{D}$, that is
\[
V^\alpha(x,y) = \inf_{T, \varphi} \{ S^\alpha_{0,T}(\varphi):
\varphi \in C([0,T], \overline{D}),  \varphi(0) = x, \varphi(T) =
y \},~~x,y \in \overline{D}.
\]

\section{Asymptotics of the Solution}
\label{asc}
\subsection{Formulation of the Result} \label{asc11}
 Recall that $(b(x), n(x)) < 0$ for $x
\in
\partial D$, where $n(x)$ is the exterior normal to the boundary
of $D$. We shall assume that there is an equilibrium  point $x_0
\in D$ for the vector field $b$, and that all the trajectories of
the dynamical system $\dot{x}(t) = b(x(t))$ starting in $D$ are
attracted to $x_0$. We also assume that there is $r >0$ such that
$(b(x), x-x_0) \leq -c|x-x_0|^2$ for some positive constant $c$
and all $x$ in the $r$-neighborhood of $x_0$.

Let $\delta > 0$, $D^\delta = \{x: x \in D,~{\rm dist}(x ,\partial
D)
> \delta \}$, and $u^\varepsilon$ be the solution of
(\ref{eq1})-(\ref{bc}). We shall be interested in the asymptotic
behavior of $u^\varepsilon(\exp({\lambda /\varepsilon^{2}}), x)$,
where $\lambda$ is fixed, $x \in D^\delta$, and~$\varepsilon
\downarrow 0$.

Let
\[
g_{\rm min} = \min_{x \in \overline{D}} g(x),~~g_{\rm max} =
\max_{x \in \overline{D}} g(x),~~ g_1 = \min_{x \in
\partial D} g(x),~~ g_2 = \max_{x
\in
\partial D} g(x).
\]
Thus $[g_1,g_2] \subseteq [g_{\rm min}, g_{\rm max}]$.  Let $M:
[g_{\rm min}, g_{\rm max}] \rightarrow \mathbb{R}$ be defined by
\begin{equation} \label{eqT}
 M(c) = \min_{x \in
\partial D} {V}^{ a(\cdot, c)}(x_0,x),
\end{equation}
where $a(x, c)$ is extended to an arbitrary bounded Lipschitz
continuous function satisfying ${k|\xi|^2 \leq \sum_{i,j =1}^d
a_{ij}(x, c) \xi_i \xi_j}$, $\xi \in \mathbb{R}^d$, $x \in
\mathbb{R}^d \setminus D$.

We next make some assumptions about the quasi-potential. It is not
difficult to see that these assumptions are satisfied by a
quasi-potential corresponding to generic~$a$ and~$b$.

 We shall assume that for all but finitely many points $c
\in [g_{\rm min}, g_{\rm max}]$ the minimum in~(\ref{eqT}) is
attained at a single point which will be denoted by $x^*(c)$. We
assume that in the remaining points $c^{1},...,c^{k}$ the minimum
is attained at two points of the boundary.
 In
this case the function  $x^*: [g_{\rm min}, g_{\rm max}]
\rightarrow
\partial D$ is piece-wise continuous and has left and right limits
at the points of discontinuity, as follows from the formula for
the quasi-potential. Let $ x_1^*(c^i) = \lim_{c \uparrow c^i}
x^*(c)$ if $c^i \neq g_{\rm min}$ and $x_2^*(c^i) = \lim_{c
\downarrow c^i} x^*(c)$ if $c^i \neq g_{\rm max}$, $1 \leq i \leq
k$. If $c^i = g_{\rm min}$, we define $x_1^*(c^i)$ as the point
distinct from $x_2^*(c^i)$ where the minimum of the
quasi-potential is attained, and similarly we define $x_2^*(c^i)$
if $c^i = g_{\rm max}$ as the point distinct from $x_1^*(c^i)$
where the minimum of the quasi-potential is attained.

We assume that $ x_1^*(c^i) \neq x_2^*(c^i)$, $1 \leq i \leq k$
(thus  $\lim_{c \uparrow c^i} x^*(c) \neq \lim_{c \downarrow c^i}
x^*(c)$ if $c^i$ is an interior point of $[g_{\rm min}, g_{\rm
max}]$). Define $G_1(c^i) = g(x_1^*(c^i))$ and $G_2(c^i) =
g(x_2^*(c^i))$.  We can now define the piece-wise continuous
function $G: [g_{\rm min}, g_{\rm max}] \rightarrow [g_1,g_2]$ via
\[
G(c) = g(x^*(c)),~~ c \in [g_{\rm min}, g_{\rm max}]\setminus
\{c^1,...,c^k\},~~~G(c^i) = G_1(c^i),~~1 \leq i \leq k.
\]
Let $c_0 = g(x_0)$ and define $c_1$ as follows:

If $G(c_0) \geq c_0$, then $c_1 = \inf \{c: c \geq c_0, G(c) \leq
c\}$.

 If $G(c_0) \leq
c_0$, then $c_1 = \sup \{c: c \leq c_0, G(c) \geq c\}$.
\\
Note that $c_1 \in [g_1,g_2]$ since $G([g_{\rm min}, g_{\rm max}])
\subseteq [g_1,g_2]$. We shall require that the graph of $G$ pass
from the left of the diagonal to the right of the diagonal at
$c_1$. More precisely, we shall assume that if $c_1 > g_{\rm
min}$, then for every $\delta_0
> 0$ there exists $\delta \in (0, \delta_0]$ such that
\[
G(c_1 - \delta) > c_1 - \delta,
\]
and if $c_1 < g_{\rm max}$, then for every $\delta_0
> 0$ there exists $\delta \in (0, \delta_0]$ such that
\[
G(c_1 + \delta) < c_1 + \delta.
\]
We also require that $c_0$ not coincide with any of the points of
discontinuity  $c^i$ for which $G_1(c^i) \leq c^i \leq G_2(c^i)$.

Let $\lambda \in (0,\infty)$ and define function $c(\lambda)$ as
follows:

For $0 < \lambda < M(c_0)$, let $c(\lambda) = c_0$.

For $\lambda \geq M(c_0)$ and $c_1 = c_0$, let $c(\lambda) = c_0$.

For $\lambda \geq M(c_0)$ and $c_1 > c_0$, let $ c(\lambda) = \min
\{c_1, \min\{c: c \in [c_0, c_1], M(c) = \lambda \} \}$.

For $\lambda \geq M(c_0)$ and $c_1 < c_0$, let $ c(\lambda) =
\max\{c_1, \max\{c: c \in [c_1, c_0], M(c) = \lambda \} \}$.

Here we use the convention that the minimum of an empty set is
$+\infty$ and the maximum of an empty set is $-\infty$. (See
Picture 1, where the thick line represents the graph of the
function $c(\lambda)$ and $\lambda'$ is a point of discontinuity
for the function $c(\lambda)$.)

We also define $\lambda_{\rm max} = \sup_{c \in [c_0,c_1]} M(c)$
if $c_1 \geq c_0$ and $\lambda_{\rm max} = \sup_{c \in [c_1,c_0]}
M(c)$ if $c_1 \leq c_0$.

\begin{figure}[htbp]
 \label{pic1}
  \begin{center}
    \begin{psfrags}
     \includegraphics[height=3.8in, width= 4.5in,angle=0]{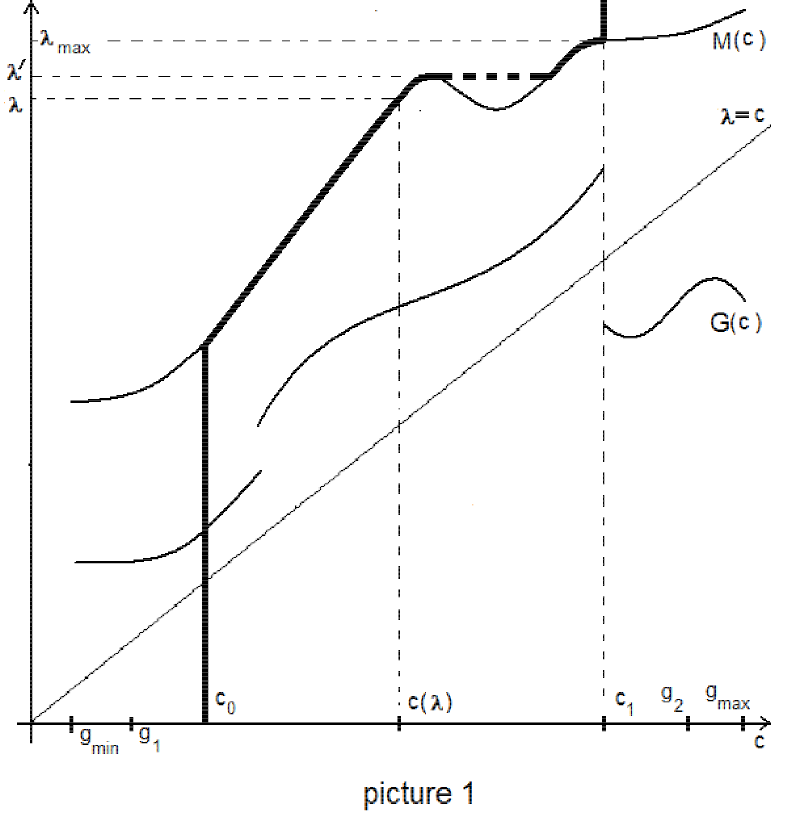}
    \end{psfrags}
    \centerline{the thick line represents the graph of $c(\lambda)$}
  \end{center}
\end{figure}

\begin{theorem} \label{mte}
Let the above assumptions concerning the differential operator
$L^\varepsilon$ and the function $G$ be satisfied. Suppose that
the function $c(\lambda)$ is continuous at a point $\lambda \in
(0, \infty)$. Then for every $\delta > 0$ the following limit
\[
\lim_{\varepsilon \downarrow 0} u^\varepsilon(\exp({\lambda
/\varepsilon^{2}}), x) = c(\lambda)
\]
is uniform in $x \in D^\delta$, where $u^\varepsilon$ is the
solution to (\ref{eq1})-(\ref{bc}).
\end{theorem}

\noindent {\bf Remark.} From Theorem~\ref{mte} and the definition
of the function $c(\lambda)$ it follows that
\[
\lim_{\varepsilon \downarrow 0} u^\varepsilon(\exp({\lambda
/\varepsilon^{2}}), x) = c_1
\]
uniformly in $x \in D^\delta$ if $\lambda > \lambda_{\rm max}$.
Moreover, from the proof of Theorem~\ref{mte} provided below it
easily follows that the limit is uniform in $(x, \lambda) \in
D^\delta \times [\overline{\lambda},\infty)$ for each
$\overline{\lambda} > \lambda_{\rm max}$. Therefore, for each
$\delta > 0$ and $\overline{\lambda} > \lambda_{\rm max}$ there is
$\varepsilon_0 > 0$ such that
\[
|u^\varepsilon(t, x) - c_1| \leq \delta
\]
whenever $\varepsilon \in (0,\varepsilon_0)$, $x \in D^\delta$ and
$t \geq \exp({\overline{\lambda} /\varepsilon^{2}})$.

It is important to note that with the boundary values of $g$
fixed, the limit $c_1$ may still depend on the initial function
through its value $c_0$ at the equilibrium point. In the generic
case, when the interval $[g_{\rm min}, g_{\rm max}]$ can be
represented as a finite union of intervals $I_1 \cup ... \cup
I_m$, such that on the interior of each of the intervals the
function $G(c) - c$ is either strictly positive or strictly
negative, the values of $c_0$ belonging to the interior of the
same interval will correspond to the same value of $c_1$.
\\

The proof of Theorem~\ref{mte} will use some properties of
diffusion processes stated in the following section.
Theorem~\ref{mte} implies various results concerning the exit
problem and metastability for the process $X^{t,x,\varepsilon}_s$
defined above. These questions will be considered in
Sections~\ref{exitprob} and \ref{genex}.

\subsection{Properties of the Diffusion Processes} \label{mmjj}

In this section we shall consider diffusion processes which are
somewhat more general than those introduced in Section~\ref{af}.
Namely, we will allow the diffusion matrix to be time dependent
but close to a matrix-valued function of the space variable (that
does not depend on time). The results stated in this section
easily follow from the arguments of \cite{FW}, Chapter~4 (see also
the companion paper \cite{FK2} for more details).

 Let $\alpha^\varepsilon$ be a symmetric $d \times d$ matrix whose elements
$\alpha^\varepsilon_{ij}$ are Lipschitz continuous on $
\mathbb{R}^+ \times \mathbb{R}^d  $ and satisfy
\begin{equation} \label{elliptic1}
k|\xi|^2 \leq \sum_{i,j =1}^d \alpha^\varepsilon_{ij}(t, x) \xi_i
\xi_j \leq K |\xi|^2,~~ (t,x) \in \mathbb{R}^+ \times
\mathbb{R}^d,~~\xi \in \mathbb{R}^d,
\end{equation}
where $k$ and $K$ are positive constants. Let $\sigma^\varepsilon$
be a square matrix such that $\alpha^\varepsilon =
\sigma^\varepsilon (\sigma^\varepsilon)^*$. We choose
$\sigma^\varepsilon$ in such a way that $\sigma^\varepsilon_{ij}$
are also bounded and Lipschitz continuous.

Let $X^{x,\varepsilon}_t$ satisfy $X^{x, \varepsilon}_0 = x$ and
\begin{equation} \label{diffpt}
d X^{x,\varepsilon}_t = b(  X^{x,\varepsilon}_t) d t +\varepsilon
\sigma^\varepsilon (t, X^{x,\varepsilon}_t) d W_t,
\end{equation}
 where $b$ is a bounded Lipschitz continuous vector field
on $ \mathbb{R}^d$ satisfying the assumptions stated in
Section~\ref{asc11}. Clearly, the law of this process depends on
$\sigma^\varepsilon$ only through $\alpha^\varepsilon =
\sigma^\varepsilon (\sigma^\varepsilon)^*$.

For $x \in D$, let $\tau^\varepsilon$ be the first time when the
process reaches the boundary of $D$. Thus $X^{
x,\varepsilon}_{\tau^\varepsilon}$ is the location of the first
exit of the process $ X^{ x,\varepsilon}_t$ from the domain $D$.
If $\alpha^\varepsilon$ is close to a function which does not
depend on time, then the asymptotics, as $\varepsilon \downarrow
0$, of $X^{ x,\varepsilon}_{\tau^\varepsilon}$ and
$\tau^\varepsilon$ can be described in terms of the
quasi-potential.

More precisely, let $\overline{\sigma}$ be a bounded Lipschitz
continuous matrix valued function on $ \mathbb{R}^d$ (with a
Lipschitz constant $L$) such that
\begin{equation} \label{kK}
k|\xi|^2 \leq \sum_{i,j =1}^d \overline{\alpha}_{ij}( x) \xi_i
\xi_j \leq K |\xi|^2,~~x \in \mathbb{R}^d,~~ \xi \in \mathbb{R}^d,
\end{equation}
 where $\overline{\alpha } =
\overline{\sigma} \overline{\sigma}^*$. Let $ \mathcal{A}$ be the
set of points in $\partial D$ at which $\min_{x
\in
\partial D} {V}^{\overline{\alpha}}(x_0,x)$ is attained. This
 minimum will be denoted by $v$.


\begin{lemma}
\label{au1} Suppose $\overline{\alpha}$ is as above with a fixed
Lipschitz constant $L$, and suppose that the positive constants
$k$ and $K$ are fixed.
 For every $\delta > 0$ there is positive $\varkappa$
  and a function $\rho: \mathbb{R}^+ \rightarrow \mathbb{R}^+$
with $\lim_{\varepsilon \downarrow 0} \rho(\varepsilon) = 0$, such
that for every
 $\alpha^\varepsilon$ that is Lipschitz continuous,
satisfies (\ref{elliptic1}) and
\begin{equation} \label{close}
 \sup_{(t,x) \in \mathbb{R}^+
\times D^\varkappa} | \alpha^\varepsilon_{ij}(t, x)-
\overline{\alpha}_{ij}(x)| \leq \varkappa,
\end{equation}
 and every $x \in D^\delta$ we have:

(A)~~~~~~~~~~~~~~~~~~~~~~ $\mathrm{P}( \tau^\varepsilon
 \leq \exp((v+\delta)/\varepsilon^2)) \geq 1 -
\rho(\varepsilon)$,

(B)~~~~~~~~~~~~~~~~~~~~~ $ \mathrm{P}(\tau^{\varepsilon} \geq
\exp((v - \delta)/\varepsilon^2)) \geq 1 - \rho(\varepsilon),$

(C)~~~~~~~~~~~~~~~~~~~~~~~~~~~$ \mathrm{P}({\rm
dist}(X^{x,\varepsilon}_{\tau^{\varepsilon}}, \mathcal{A}) \leq
\delta) \geq 1 - \rho(\varepsilon)$.
%
\end{lemma}
The next lemma only requires the boundedness of the quadratic form
$\alpha^\varepsilon$ from above and below.
\begin{lemma} \label{sec2}
Suppose that positive constants $k$ and $K$ are fixed. There
exists $v_0 > 0$ such that for every $0 < \delta < v_0$ there is a
function $\rho: \mathbb{R}^+ \rightarrow \mathbb{R}^+$ with
$\lim_{\varepsilon \downarrow 0} \rho(\varepsilon) = 0$, such that
for every
 $\alpha^\varepsilon$ that is Lipschitz continuous and
satisfies (\ref{elliptic1}) and every $x \in D^\delta$ we have:

(A) ~~~~~~~~~~~~~~~~~~~~~~~$ \mathrm{P}(\tau^\varepsilon \geq
\exp(v_0/\varepsilon^2)) \geq 1 - \rho(\varepsilon),$

(B) ~~~~~~~~~~~~~~~~$\inf_{t \in  [\exp(\delta/\varepsilon^2),
\exp(v_0/\varepsilon^2)]} \mathrm{P}(|X^{x,\varepsilon}_t - x_0|
\leq \delta) \geq 1 - \rho(\varepsilon),$

(C)~~~~~~~~~~~~~~~~$ \inf_{t \in [0, \exp(v_0/\varepsilon^2)]}
\mathrm{P}(|X^{ x_0,\varepsilon}_t - x_0| \leq \delta) \geq 1 -
\rho(\varepsilon)$.
\end{lemma}
An easy corollary of this lemma is that at an exponential time the
process either can be found in a small neighborhood of $x_0$ or
has earlier crossed the boundary of the domain.
\begin{corollary} \label{co1}
Suppose that positive constants $k$ and $K$ are fixed. For every
$\delta > 0$ there  is a function $\rho: \mathbb{R}^+ \rightarrow
\mathbb{R}^+$  with $\lim_{\varepsilon \downarrow 0}
\rho(\varepsilon) = 0$, such that  for every
 $\alpha^\varepsilon$ that is Lipschitz continuous and
satisfies (\ref{elliptic1}), every $x \in D$ and $t \geq
\exp(\delta/\varepsilon^2)$ we have:
\[
\mathrm{P}(|X^{x,\varepsilon}_t - x_0| \leq \delta~~{\rm
or}~~\tau^\varepsilon \leq t) \geq 1 - \rho(\varepsilon).
\]
\end{corollary}
\proof Let $\delta_1 > 0$ be sufficiently small so that there is a
domain $\widetilde{D}$ with smooth boundary  such that
$\widetilde{D}^{\delta_1} = D$. If the process does not reach
$\partial D$ by the time $t - \exp(\delta/\varepsilon^2)$, then we
can apply Part (B) of Lemma~\ref{sec2} to the domain
$\widetilde{D}$ and the process starting at
$X^{x,\varepsilon}_{t-\exp(\delta/\varepsilon^2)}$, and the result
follows from the Markov property. \qed

\subsection{Preliminary Lemmas}
The next step in the proof of Theorem~\ref{mte} is to establish
that $u^\varepsilon(\exp({t /\varepsilon^{2}}), x)$ is nearly
constant on $D^\delta$ if $t > 0$ is fixed and $\varepsilon$ is
sufficiently small. This is accomplished in Lemma~\ref{letwo}
below.
\begin{lemma} \label{small}
For every positive $t_0$ and $\delta$  there are positive $R$ and
$\varepsilon_0$ such that
\begin{equation} \label{epsi}
 |u^\varepsilon(t, x) -
u^\varepsilon(t,x_0)| \leq \delta
\end{equation}
whenever $|x - x_0| \leq R \varepsilon$,   $\varepsilon \leq
\varepsilon_0$ and $t \geq t_0$.
\end{lemma}
\proof
Let $v^\varepsilon(t, y) = u^\varepsilon(t, x_0+\varepsilon y)$,
$t \in (0,\infty)$, $|y| \leq 2$. Let $Q = (0,\infty) \times
B_{2}(0)$ and $Q_0 = (t_0, \infty) \times B_1(0)$, where $B_r(0)$
is the ball of radius $r$ centered at the origin. Then
$v^\varepsilon$ satisfies the following partial differential
equation:
\[
\frac{\partial v^\varepsilon(t,y)}{\partial t}  = \frac{1}{2}
\sum_{i,j=1}^d \widetilde{a}_{ij}(t, y) \frac{
\partial^2 v^\varepsilon (t,y)}{\partial y_i \partial y_j} + \frac{b(x_0 + \varepsilon
y)}{\varepsilon} \cdot \nabla_y v^\varepsilon (t,y),~~(t,y) \in Q.
\]
Here $\widetilde{a}_{ij}(t,y) = a_{ij} (x_0 + \varepsilon y,
u^\varepsilon(t, x_0 + \varepsilon y))$ are uniformly bounded in
$\varepsilon$ and satisfy $k|\xi|^2 \leq \sum_{i,j =1}^d
\widetilde{a}_{ij}(t,y) \xi_i \xi_j$, $t \in (0,\infty)$, $y \in
B_{2}(0)$, $\xi \in \mathbb{R}^d$. Moreover, $\sup_{y \in
B_{2}(0)} | {b(x_0 + \varepsilon y)}/{\varepsilon}|$ is bounded
uniformly in $\varepsilon$ and $ |\nabla_y
\widetilde{a}_{ij}(t,y)|$ can be estimated from above by a
constant times $1+|\nabla_y v^\varepsilon (t,y)|$ for $(t,y) \in
Q$, uniformly in $\varepsilon$. Since the distance between $Q_0$
and the boundary of $Q$ is positive and $|v^\varepsilon|$ is
uniformly bounded in $Q$ by $\max(|g_1|, |g_2|)$, we can apply the
a priori estimate (see Theorem 4, Chapter 5.2 of \cite{K} or
Theorem 6, Chapter~6.2 of \cite{K}) to bound $\sup_{(t,x) \in Q_0}
|\nabla_y v^\varepsilon|$ by  a constant $C$ independent of
$\varepsilon$. This implies that $|u^\varepsilon(t,
x_0+\varepsilon y)|$ is bounded by $C R$ when $|y| \leq R$ if $R$
is a positive constant. It remains to take $R$ such that $C R \leq
\delta$. \qed

We'll need the following simple lemma about diffusion processes
with the drift directed towards the origin.
\begin{lemma} \label{dipro}
Let $\widetilde{b}$ be a $ {C}^2$ smooth vector field on $
\mathbb{R}^d$ such that $(\widetilde{b}(x),x) \leq -k_1(x,x)$ for
some positive $k_1$ and all $x \in \mathbb{R}^d$. Let
${\sigma^\varepsilon}(t,x)$ be a Lipschitz continuous function
such that ${\alpha^\varepsilon} = {\sigma^\varepsilon}
{(\sigma^\varepsilon)}^*$ satisfies (\ref{elliptic1}).
Let $Y^{x,\varepsilon}_t$ be the process starting at $x$ that
satisfies
\[
d Y^{x,\varepsilon}_t = \widetilde{b}( Y^{x,\varepsilon}_t ) d t+
\varepsilon {\sigma^\varepsilon}(t, Y^{x,\varepsilon}_t) d W_t.
\]
Then for every $r, R > 0$ there are positive $\gamma$, $s_0$, and
$\varepsilon_0$, which depend on $\widetilde{b}$ and
${\sigma^\varepsilon}$ only through $k_1$, $k$, and $K$ such that
 \begin{equation} \label{inside2}
 \mathrm{P}(Y^{x,\varepsilon}_{s_0 |\ln \varepsilon|} \in B_{R
\varepsilon}(0)) \geq \gamma
\end{equation}
 holds for $x \in B_r(0)$ and
$0< \varepsilon \leq \varepsilon_0$.
\end{lemma}
\proof First, let us show that the probability that the process
enters a larger ball in time $s_0 |\ln \varepsilon| - 1$ is
bounded from below. Let $h: \mathbb{R} \rightarrow [0,1]$ be a
smooth even function with negative derivative on $(1/2,1)$, such
that $h(x) = 1$ for $0 \leq x \leq 1/2$  and $h(x) = 0$ for $x
\geq 1$. Let
\[
f(t,x) = h(|x|\exp(-2 k_1 t)/R_1 \varepsilon),
\]
where $R_1 > 0$ will be specified below. If $R_1$ and $s_0$ are
sufficiently large and $\varepsilon_0$ is sufficiently small, then
$f(s_0 |\ln \varepsilon| -1 ,x) =1$ for $x \in B_r(0)$, $0 <
\varepsilon \leq \varepsilon_0$. By the Ito formula,
\[
\mathrm{P}(Y^{x,\varepsilon}_{s_0 |\ln \varepsilon|-1} \in B_{R_1
\varepsilon}(0)) \geq \mathrm{E}f(0, Y^{x,\varepsilon}_{s_0 |\ln
\varepsilon|-1})=
\]
\[
\mathrm{E}f(s_0 |\ln \varepsilon|-1, x) + \mathrm{E} \int_0^{ s_0
|\ln \varepsilon|-1} (\mathcal{L}^\varepsilon f - \frac{\partial
f}{\partial t})(s_0 |\ln \varepsilon|-1 -s, Y^{x,\varepsilon}_s) d
s,
\]
where $ \mathcal{L}^\varepsilon$ is the generator of the process
$Y^{x,\varepsilon}_t$. In order to estimate the integral in the
right hand side, we note that
\[
 (\mathcal{L}^\varepsilon f - \frac{\partial f}{\partial
t}) (t,x) \geq - C \max_{x \in [0,1]}(h''(x)) \exp(-2k_1 t)/R_1^2,
\]
where the constant $C$ depends on $k_1$ and $K$. By taking $R_1$
sufficiently large, we can bound the expectation of the integral
from below by $-1/2$. We have thus demonstrated that
\[
 \mathrm{P}(Y^{x,\varepsilon}_{s_0 |\ln \varepsilon| -1} \in B_{R_1
\varepsilon}(0)) \geq 1/2
\]
holds for $x \in B_r(0)$ and $0< \varepsilon \leq \varepsilon_0$.
This, together with the Markov property, will imply
(\ref{inside2}) if we show that there is $\gamma > 0$ such that
\[
 \mathrm{P}(Y^{x,\varepsilon}_{1} \in B_{R
\varepsilon}(0)) \geq 2 \gamma
\]
holds for $x \in  B_{R_1 \varepsilon}(0)$ and $0< \varepsilon \leq
\varepsilon_0$. The latter inequality is a consequence of the
Aronson estimate (\cite{Ar}).  Indeed, we can make the same change
the variables in the generator of the process as in the proof of
Lemma~\ref{small}, thus obtaining a uniformly elliptic operator,
and therefore the Aronson estimate is applicable. \qed
\begin{lemma} \label{tsev}
Let the process $X^{x,\varepsilon}_t$ satisfy (\ref{diffpt}) with
the coefficients satisfying the assumptions stated in Section
\ref{mmjj}. Suppose that positive constants $k$, $K$, and $R$ are
fixed. There exist $\gamma > 0$ and  $v_0 > 0$ such that for every
$0 < \delta < v_0$ there is $\varepsilon_0 > 0$, such that for
every
 $\alpha^\varepsilon$ that is Lipschitz continuous and
satisfies (\ref{elliptic1}) and every $x \in D^\delta$ we have:
\begin{equation} \label{chii}
\inf_{t \in  [\exp(\delta/\varepsilon^2),
\exp(v_0/\varepsilon^2)]} \mathrm{P}(| X^{x,\varepsilon}_t - x_0|
\leq R \varepsilon) \geq \gamma
\end{equation}
when $\varepsilon \leq \varepsilon_0$.
\end{lemma}
\proof This lemma easily follows from Lemmas \ref{sec2},
\ref{dipro}, and the Markov property of the process. \qed
\\
\\
Finally, we show that the function $u^\varepsilon$ is nearly
constant on the domain for each fixed, sufficiently large, value
of time.
\begin{lemma} \label{letwo}
For every positive $\lambda_0$ and $\delta$ there is positive
$\varepsilon_0$ such that
\begin{equation} \label{kap} |u^\varepsilon(\exp({\lambda
/\varepsilon^{2}}), x) - u^\varepsilon(\exp({\lambda
/\varepsilon^{2}}),x_0)| \leq \delta
\end{equation}
whenever $x \in D^\delta$, $\varepsilon \leq \varepsilon_0$ and
$\lambda \geq \lambda_0$.
\end{lemma}
\proof Take $t_0 = 1$ and find $\varepsilon_0$ and $R$ such that
(\ref{epsi}) holds with $\delta/4$ instead of $\delta$ in the
right hand side, provided that $|x - x_0| \leq R \varepsilon$,
$\varepsilon \leq \varepsilon_0$ and $t \geq t_0$. By
Lemma~\ref{tsev} we can take $\gamma > 0$ and  $v_0 < \lambda_0$
such that for each $\delta' > 0$ we have
\[
\mathrm{P}(| X^{x,\varepsilon}_{ \exp(v_0/\varepsilon^2)} - x_0|
\leq R \varepsilon) \geq \gamma
\]
for all $x \in D^{\delta'}$ and all sufficiently small
$\varepsilon$. Observe that
\[
\sup_{t \geq 0,~x,y \in D} |u^\varepsilon(t,x) -
u^\varepsilon(t,y) | \leq \sup_{x,y \in D} |g(x) - g(y)|.
\]
Choose an integer
\begin{equation} \label{iin}
N \geq 1 + 6 \sup_{x,y \in D} |g(x) - g(y)| / (\delta \gamma)
\end{equation}
and consider the sequence of domains $D^\delta \subset
D^{\delta/2} \subset...\subset D^{\delta/N}$. Let
\[
q_m = \sup_{x,y \in D^{\delta/m}} |u^\varepsilon(e^{\frac{\lambda}
{\varepsilon^{2}}} - (m-1) e^{\frac{v_0}{\varepsilon^{2}}},x) -
u^\varepsilon(e^{\frac{\lambda}{\varepsilon^{2}}} - (m-1)
e^{\frac{v_0}{\varepsilon^{2}}},y)|.
\]
We claim that
\begin{equation} \label{recur}
q_m \leq  (1 - \frac{\gamma}{2}) q_{m+1} + \frac{\gamma
\delta}{3},~~~~m=1,...,N.
\end{equation}
for all sufficiently small $\varepsilon$. Since $q_N \leq
\sup_{x,y \in D} |g(x) - g(y)|$, the inequalities (\ref{iin}) and
(\ref{recur}) imply that $q_1 \leq \delta$, which gives
(\ref{kap}). It remains to prove (\ref{recur}).

We suppress the dependence on $m$ and $\varepsilon$ in the
notation for the process $Y^x_s = Y^{m,x,\varepsilon}_s$, $s \in
[0,\exp({v_0 /\varepsilon^{2}})]$, that starts at $x \in
D^{\delta/m}$ and satisfies
\begin{equation} \label{pro}
d Y^{x}_s = b( Y^{x}_s) d s+ \varepsilon \sigma(Y^{x}_s,
u^\varepsilon(e^{\frac{\lambda} {\varepsilon^{2}}} - (m-1)
e^{\frac{v_0}{\varepsilon^{2}}}-s, Y^{x}_s)) d W_s,~~s \leq \tau^x
\wedge \exp({v_0 /\varepsilon^{2}}),
\end{equation}
\[
\tau^x = \min\{s: Y^{x}_s \in
\partial D^{\delta/(m+1)}\},~~~~~ Y^{x}_s =
 Y^{x}_{\tau^x},~~\tau^x \leq s
 \leq \exp({v_0 /\varepsilon^{2}}),
\]
where $\sigma_{ij} \in C^2 (D \times \mathbb{R})$, $1 \leq i,j
\leq d$, are Lipschitz continuous and such that $\sigma \sigma^* =
a$. (If the minimum in the definition of $\tau^x$ is taken over an
empty set, then it is considered to be equal to $+ \infty$.)

Let $\Omega$ denote the probability space on which the diffusion
process is defined, and consider its partition into the following
disjoint events:
\[
G_1^x = \{ | Y^{x}_{ \exp(v_0/\varepsilon^2)} - x_0| \leq R
\varepsilon, ~~ \tau^x \geq \exp({v_0 /\varepsilon^{2}}) \},
\]
\[
G_2^x = \{ | Y^{x}_{ \exp(v_0/\varepsilon^2)} - x_0| > R
\varepsilon, ~~ \tau^x \geq \exp({v_0 /\varepsilon^{2}}) \},
\]
\[
G_3^x = \{  \tau^x < \exp({v_0 /\varepsilon^{2}}) \}.
\]
From our construction and  Lemma \ref{sec2} it follows that $
\mathrm{P}(G_1^x) \geq \gamma/2$ for all sufficiently
small~$\varepsilon$. Let us take an arbitrary event $ F_1^x
\subseteq G_1^x$ such that $ \mathrm{P}(F_1^x) = \gamma/2$. From
Lemma~\ref{sec2} it follows that $ \mathrm{P}(G_3^x) \leq \gamma
\delta/ (1+24\sup_{x,y \in D} |g(x) - g(y)| )$ for all
sufficiently small $\varepsilon$. We can assume that $\gamma$ is
sufficiently small so that $ \gamma \delta/ (1+24\sup_{x,y \in D}
|g(x) - g(y)| ) + \gamma/2 < 1$. Let us take an arbitrary event
$F_3^x$ such that $\mathrm{P}(F_3^x) = \gamma \delta/ (1+24
\sup_{x,y \in D} |g(x) - g(y)| )$ and $ G_3^x \subseteq F_3^x
\subseteq \Omega \setminus F_1^x$. Let $F_2^x = \Omega \setminus
(F_1^x \cup F_3^x)$.

By the Feynman-Kac formula,
\[
u^\varepsilon(e^{\frac{\lambda} {\varepsilon^{2}}} - (m-1)
e^{\frac{v_0}{\varepsilon^{2}}},x)  =
\mathrm{E}u^\varepsilon(e^{\frac{\lambda} {\varepsilon^{2}}} -
(m-1)e^{\frac{v_0}{\varepsilon^{2}}}-\tau^x,Y^{x}_{\tau^x})=
\]
\[
\mathrm{E}[u^\varepsilon(e^{\frac{\lambda} {\varepsilon^{2}}} -
(m-1)e^{\frac{v_0}{\varepsilon^{2}}}-\tau^x,Y^{x}_{\tau^x})
\chi_{F^x_1}]  + \mathrm{E}[u^\varepsilon(e^{\frac{\lambda}
{\varepsilon^{2}}} -
(m-1)e^{\frac{v_0}{\varepsilon^{2}}}-\tau^x,Y^{x}_{\tau^x})
\chi_{F^x_2}] +
\]
\[
\mathrm{E}[u^\varepsilon(e^{\frac{\lambda} {\varepsilon^{2}}} -
(m-1)e^{\frac{v_0}{\varepsilon^{2}}}-\tau^x,Y^{x}_{\tau^x})
\chi_{F^x_3}] = I^x_1 + I^x_2 + I^x_3.
\]
Therefore,
\[
q_m \leq \sup_{x,y \in D_m} (|I^x_1  - I^y_1 |+ |I^x_2  - I^y_2| +
|I^x_3  - I^y_3|).
\]
The term $|I^x_1  - I^y_1 |$ is estimated from above by $\gamma
\delta /4$ since the value of the function $u^\varepsilon$ inside
the $R \varepsilon$ neighborhood of $x_0$ doesn't vary by more
than $\delta/2$. The term $  |I^x_2  - I^y_2| $ can be estimated
from above by $(1- \gamma/2) q_{m+1}$ since $\mathrm{P}(F^x_2) =
\mathrm{P}(F^y_2) \leq 1 - \gamma/2$. Similarly, the term $ |I^x_3
- I^y_3|$ is estimated from above by $\gamma \delta /12$.
Combining the terms, yields (\ref{recur}), which completes the
proof of the lemma. \qed

\subsection{Proof of the Theorem on the Asymptotics of the Solution}
In this section we prove Theorem~\ref{mte}. First, we examine the
behavior of $u^\varepsilon$ for times which are small in the
logarithmic scale.
\begin{lemma} \label{simp}
There is a positive $v_0$ such that for every $0 < \delta < v_0$
there is  $\varepsilon_0 > 0$ such that
\[
 |u^\varepsilon(\exp({\lambda
/\varepsilon^{2}}), x) - g(x_0)| \leq \delta
\]
whenever $x \in D^\delta$, $0 < \varepsilon \leq \varepsilon_0$
and $\delta \leq \lambda \leq v_0$.
\end{lemma}
\proof This lemma immediately follows from Lemma~\ref{sec2}. \qed
\\

The next three lemmas, central to the proof of Theorem~\ref{mte},
rule out certain types of behavior for the function
$u^\varepsilon$.
\begin{lemma} \label{fo}
Suppose that $a_1 \leq \mu_n < \lambda_n \leq a_2$ for some
constants $a_1, a_2 >0$, $\varepsilon_n \downarrow 0$ as $n
\rightarrow \infty$, and
\[
u^{\varepsilon_n}(\exp(\mu_n/\varepsilon_n^2), x_0) = \beta_1,~~
u^{\varepsilon_n}(\exp(\lambda_n/\varepsilon_n^2), x_0) = \beta_2
\]
with $\beta_1 \neq \beta_2$. Then there is $\delta > 0$ such that
\begin{equation} \label{gap}
\exp(\lambda_n/\varepsilon_n^2) - \exp(\mu_n/\varepsilon_n^2) \geq
\exp(\delta/\varepsilon_n^2)
\end{equation}
for all large enough $n$.
\end{lemma}
\proof Consider the process $ X^{\lambda_n,x_0,\varepsilon_n}_s$
given by~(\ref{pro}) with $\tau^{\varepsilon_n}$ being the first
time when this process reaches the boundary of $D$. Define
\[
\overline{\tau}^{\varepsilon_n} = \min(\tau^{{\varepsilon_n}},
\exp(\lambda_n/\varepsilon_n^2)-\exp(\mu_n/\varepsilon_n^2)).
\]
 Then
\[
u^{\varepsilon_n}(\exp({\lambda_n /\varepsilon_n^{2}}),x_0) =
\mathrm{E} u^{\varepsilon_n}(\exp({\lambda_n
/\varepsilon_n^{2}})-\overline{\tau}^{\varepsilon_n},
X^{\lambda_n,x_0,\varepsilon_n}_{\overline{\tau}^{\varepsilon_n}}).
\]
The left hand side in this formula is equal to $\beta_2$. If
(\ref{gap}) does not hold, then by Part~(C) of Lemma~\ref{sec2}
and Lemma~\ref{letwo} the right hand side can be made arbitrarily
close to $\beta_1$ along a subsequence. \qed

\begin{lemma} \label{tet}
Suppose that $a_1 \leq \mu_n < \lambda_n \leq a_2$ for some
constants $a_1, a_2 >0$, $\varepsilon_n \downarrow 0$ as $n
\rightarrow \infty$, and
\begin{equation} \label{ml}
u^{\varepsilon_n}(\exp(\mu_n/\varepsilon_n^2), x_0) = \beta_1,~~
u^{\varepsilon_n}(\exp(\lambda_n/\varepsilon_n^2), x_0) = \beta_2,
\end{equation}
If $g_{\rm min} < \beta_1 < \beta_2 < g_{\rm max}$, then neither
of the following is possible:

  (A) There is $\delta > 0$ such that $\lambda_n < M(\beta_2) -
  \delta$,

  (B) There is $\delta > 0$ such that $G(c) < \beta_2 - \delta$ for $c \in
  [\beta_2 -\delta, \beta_2+\delta]$.
\\
 If  $g_{\rm min} < \beta_2 < \beta_1 < g_{\rm max}$, then
neither of the following is possible:

  (A$'$) There is $\delta > 0$ such that $\lambda_n < M(\beta_2) -
  \delta$,

  (B$'$) There is $\delta > 0$ such that $G(c) > \beta_2 + \delta$ for $c \in
  [\beta_2 -\delta, \beta_2+\delta]$.
\end{lemma}
\proof  (A) Let us assume that $\lambda_n < M(\beta_2) -
  \delta$. We can find $\lambda'_n \in [\mu_n, \lambda_n]$ such that
$u^{\varepsilon_n}(\exp(\lambda'_n/\varepsilon_n^2), x_0) =
\beta_2$ and
 $u(t,x_0) \leq \beta_2$  for $t \in
[\exp(\mu_n/\varepsilon_n^2),\exp(\lambda'_n/\varepsilon_n^2)]$.

   Let $\overline{\alpha}_{ij}(x) =
a_{ij}(x, \beta_2)$. Let $\tau^{\varepsilon}$ be the first time
when the process defined in~(\ref{diffpt}) reaches the boundary of
$D$. By part (B) of Lemma~\ref{au1}, we can choose $\varkappa
> 0$ such that whenever
 $\alpha$ is Lipschitz continuous and
satisfies (\ref{elliptic1}) and (\ref{close}), we have
\begin{equation} \label{er1}
\mathrm{P}(\tau^{\varepsilon_n} \geq
\exp(\lambda'_n/\varepsilon_n^2)) \geq
\mathrm{P}(\tau^{\varepsilon_n} \geq \exp((M(\beta_2) -
\delta)/\varepsilon_n^2)) \geq 1 - \rho(\varepsilon_n),
\end{equation}
for $x \in D^\delta$, where $\rho$ does not depend on $\alpha$ and
satisfies $\lim_{\varepsilon \downarrow 0} \rho(\varepsilon) = 0$.
Choose $\beta' > 0$ such that
\[
|a_{ij}(x,\beta_2) - a_{ij}(x,\beta)| < \varkappa
\]
whenever $\beta \in [\beta_2 - 2\beta', \beta_2 + 2 \beta']$, $x
\in D^\varkappa$. Choose a sequence $\mu'_n \in [\mu_n,
\lambda'_n]$ such that $u^{\varepsilon_n}(\mu'_n/\varepsilon_n^2,
x_0) = \beta_2 - \beta'$ and $u^{\varepsilon_n}(t,x_0) \in
[\beta_2 - \beta', \beta_2]$ for $t \in
[\exp(\mu'_n/\varepsilon_n^2),\exp(\lambda'_n/\varepsilon_n^2)]$.
 By Lemma~\ref{letwo}, we have
\begin{equation} \label{as}
|a_{ij}(x,\beta_2) - a_{ij}(x,u^{\varepsilon_n}(t,x))| < \varkappa
\end{equation}
for $x \in D^\varkappa$, $t \in [\exp(\mu'_n/\varepsilon_n^2),
\exp(\lambda'_n/\varepsilon_n^2)]$, if $\varepsilon_n$ is
sufficiently small. Consider the process $
X^{\lambda'_n,x_0,\varepsilon_n}_s$ given by~(\ref{pro}), with
$\tau^{\varepsilon_n}$ now being the first time when this process
reaches the boundary of $D$, and
\begin{equation} \label{tau}
\overline{\tau}^{\varepsilon_n} = \min(\tau^{{\varepsilon_n}},
\exp(\lambda'_n/\varepsilon_n^2)-\exp(\mu'_n/\varepsilon_n^2)).
\end{equation}
 Then
\begin{equation} \label{interv}
u^{\varepsilon_n}(\exp({\lambda'_n /\varepsilon_n^{2}}),x_0) =
\mathrm{E} u^{\varepsilon_n}(\exp({\lambda'_n
/\varepsilon_n^{2}})-\overline{\tau}^{\varepsilon_n},
X^{\lambda'_n,x_0,\varepsilon_n}_{\overline{\tau}^{\varepsilon_n}}).
\end{equation}
The left hand side in this formula is equal to $\beta_2$, while
the right hand side can be made arbitrarily close to $\beta_2 -
\beta'$ by considering sufficiently small $\varepsilon_n$ due to
(\ref{er1}), Corollary~\ref{co1} (which applies due to
Lemma~\ref{fo}) and Lemma~\ref{letwo}. This leads to a
contradiction.

(B) Assume that $G(c) < \beta_2 - \delta$ for $c \in
  [\beta_2 -\delta, \beta_2+\delta]$. Let $ \mathcal{A}$ be the
  set (consisting of either one or two points) where the minimum
  of ${V}^{a(\cdot, \beta_2)}(x_0,x)$ is attained. From the definition
  of $G$ it follows that $g( \mathcal{A}) \subset (-\infty ,
  \beta_2 - \delta)$.  By part (C) of Lemma~\ref{au1}, we can
choose $\varkappa
> 0$ such that whenever
 $\alpha$ is Lipschitz continuous and
satisfies (\ref{elliptic1}) and (\ref{close}), we have
\begin{equation} \label{pp}
 \mathrm{P}(g(X^{\sigma, x,\varepsilon}_{\tau^{\varepsilon}})
\leq \beta_2-\delta/2) \geq 1 - \rho(\varepsilon)
\end{equation}
for $x \in D^\delta$ and all sufficiently small $\varepsilon$.

  As in case (A),
we can find $\beta' >0$ and $\mu_n \leq \mu'_n < \lambda'_n \leq
\lambda_n$ such that
\[
u^{\varepsilon_n}(\exp(\mu'_n/\varepsilon_n^2), x_0) = \beta_2 -
\beta',~~ u^{\varepsilon_n}(\exp(\lambda'_n/\varepsilon_n^2), x_0)
= \beta_2,
\]
and (\ref{as}) holds for $x \in D^\varkappa$, $t \in
[\exp(\mu'_n/\varepsilon_n^2), \exp(\lambda'_n/\varepsilon_n^2)]$,
if $\varepsilon_n$ is sufficiently small. We can again employ
formula (\ref{interv}) in which the left hand side is equal to
$\beta_2$. The right hand side can be written as
\[
\mathrm{E} u^{\varepsilon_n}(\exp({\lambda'_n
/\varepsilon_n^{2}})-\overline{\tau}^{\varepsilon_n},
X^{\lambda'_n,x_0,\varepsilon_n}_{\overline{\tau}^{\varepsilon_n}})=
\]
\[
\mathrm{E}(\chi_{\{\overline{\tau}^{\varepsilon_n} <
{\tau}^{\varepsilon_n} \}}
 u^{\varepsilon_n}(\exp({\mu'_n
/\varepsilon_n^{2}}),
X^{\lambda'_n,x_0,\varepsilon_n}_{\overline{\tau}^{\varepsilon_n}}))+
\mathrm{E}(\chi_{\{\overline{\tau}^{\varepsilon_n} =
{\tau}^{\varepsilon_n} \}}
g(X^{\lambda'_n,x_0,\varepsilon_n}_{{\tau}^{\varepsilon_n}})).
\]
The first term in the right hand side here can be made arbitrarily
close to $ \mathrm{P}(\overline{\tau}^{\varepsilon_n} <
{\tau}^{\varepsilon_n}) (\beta_2-\beta')$ by Corollary~\ref{co1}
and Lemma~\ref{letwo}. The second term on the right hand side can
be estimated from above for large $n$ by $
\mathrm{P}(\overline{\tau}^{\varepsilon_n} =
{\tau}^{\varepsilon_n})(\beta_2 - \delta/4)$ due to (\ref{pp}).
This leads to a contradiction.

The proof of (A$'$) and (B$'$) is completely similar to the proof
of (A) and (B). \qed

\begin{lemma} \label{rone}
Suppose that $ \lambda > 0$. If $c_1 > c_0$, then
\begin{equation} \label{fsta}
\liminf_{\varepsilon \downarrow 0}
(u^{\varepsilon}(\exp(\lambda/\varepsilon^2), x_0) -c(\lambda))
\geq 0.
\end{equation}
 If
$c_1 < c_0$, then
\[
\limsup_{\varepsilon \downarrow 0}
(u^{\varepsilon}(\exp(\lambda/\varepsilon^2), x_0) -c(\lambda))
\leq 0.
\]
\end{lemma}
\proof We shall only consider the first statement since the second
one is completely similar. Note that
\begin{equation} \label{spll}
\liminf_{\varepsilon \downarrow 0}
u^{\varepsilon}(\exp(\lambda/\varepsilon^2), x_0) \geq c_0.
\end{equation}
Indeed, otherwise by Lemma~\ref{simp} there are $\beta_2 < \beta_1
< c_0$ and sequences $v_0 \leq \mu_n < \lambda_n \leq \lambda$ and
$\varepsilon_n \downarrow 0$ such that (\ref{ml}) holds. Note that
the graph of $G$ goes above the diagonal in a neighborhood of
$c_0$, while $\beta_1$ and $\beta_2$ can be taken arbitrarily
close to $c_0$. Therefore, Part~(B$'$) of Lemma~\ref{tet} leads to
a contradiction.

Thus, if (\ref{fsta}) does not hold, then is there are $\delta >
0$ and a sequence $\varepsilon_n \downarrow 0$ such that
\[
c_0 - \delta < u^{\varepsilon_n}(\exp(\lambda/\varepsilon_n^2),
x_0) < c(\lambda)-\delta.
\]
We choose $\delta$ sufficiently small so that the graph of $G$
goes above the diagonal on the interval $[c_0 - 2\delta,
c(\lambda) - \delta/2]$.
 Take $\delta' > 0$ which
will be specified later.

For each $c \in [c_0-\delta, c(\lambda) - \delta]$, by part (A) of
Lemma~\ref{au1}, we can choose $\varkappa(c) > 0$ such that
whenever $\alpha$ is Lipschitz continuous and satisfies
(\ref{elliptic1}) and (\ref{close}) with
$\overline{\alpha}_{ij}(x) = a_{ij}(x,c)$  we have
\begin{equation} \label{hn}
\mathrm{P}( \tau^\varepsilon
 \leq \exp((M(c)+\delta')/\varepsilon^2)) \geq 1 -
\rho(\varepsilon).
\end{equation}
\begin{equation} \label{hn2}
 \mathrm{P}(g(X^{x_0,\varepsilon}_{\tau^{\varepsilon}}) \geq
\inf_{\overline{c} \in [c-\delta', c+ \delta']}G(\overline{c}) -
\delta' ) \geq 1 - \rho(\varepsilon)
\end{equation}
For  each $c \in [c_0-\delta, c(\lambda) - \delta]$, find $l(c) <
\delta'$ such that $|a_{ij}(x,c) - a_{ij}(x,\overline{c})| <
\varkappa(c)$ whenever $\overline{c} \in [c - l(c), c + l(c)]$, $x
\in D^{\varkappa(c)}$.

Choose a finite subcovering of the interval $[c_0-\delta,
c(\lambda) - \delta]$ by the intervals $(c_m - l(c_m)/2,c_m +
l(c_m)/2)$ and take $l = \min_m(l(c_m))$. Let $c_0-\delta =
\beta_0 < \beta_1<...<\beta_k = c(\lambda) - \delta$ be such that
$\beta_{i} - \beta_{i-1} \leq l/10$, $1 \leq i \leq k$.

We claim that if $i \geq 1$, $0 < \lambda' < \lambda$ and
$u^{\varepsilon_n}(\exp(\lambda/\varepsilon_n^2), x_0) \in
[\beta_{i-1}, \beta_i]$ along a subsequence, then
$u^{\varepsilon_n}(\exp(\lambda'/\varepsilon_n^2), x_0) \leq
\beta_{i-1}$ for large enough $n$ along the same subsequence. If
this were not the case, then we would have
\[
u^{\varepsilon_n}(\exp(\lambda/\varepsilon_n^2), x_0) \in
[\beta_{i-1}, \beta_i]~~{\rm
and}~~u^{\varepsilon_n}(\exp(\lambda'/\varepsilon_n^2), x_0) \geq
\beta_{i-1}
\]
along a further subsequence. Note that the function
$u^{\varepsilon_n}(t, x_0)$ must take values in the interval
$[\beta_{i-1} - l/10, \beta_i + l/10]$ for $t \in
[\exp(\lambda'/\varepsilon_n^2),\exp(\lambda/\varepsilon_n^2)]$,
otherwise Part (B$'$) of Lemma~\ref{tet} leads to a contradiction.

 By the construction above and Lemma~\ref{letwo},
 there is $m$ such that $|c_{m} - \beta_i| < \delta'$ and
\[
|a_{ij}(x,c_{m}) - a_{ij}(x,u^{\varepsilon_n}(t,x))| <
\varkappa(c_{m})
\]
for $x \in D^{\varkappa(c_{m})}$, $t \in
[\exp(\lambda'/\varepsilon_n^2), \exp(\lambda/\varepsilon_n^2)]$,
if $\varepsilon_n$ is sufficiently small.

 Consider the process $
X^{\lambda,x_0,\varepsilon_n}_s$ given by~(\ref{pro}), with
$\tau^{\varepsilon_n}$ being the first time when this process
reaches the boundary of $D$, and
\[
\overline{\tau}^{\varepsilon_n} = \min(\tau^{{\varepsilon_n}},
\exp(\lambda/\varepsilon_n^2)-\exp(\lambda'_n/\varepsilon_n^2)).
\]
 Then
\[ u^{\varepsilon_n}(\exp({\lambda /\varepsilon_n^{2}}),x_0) =
\mathrm{E} u^{\varepsilon_n}(\exp({\lambda
/\varepsilon_n^{2}})-\overline{\tau}^{\varepsilon_n},
X^{\lambda,x_0,\varepsilon_n}_{\overline{\tau}^{\varepsilon_n}}).
\]
The left hand side does not exceed $\beta_i$. If $(
\exp(\lambda/\varepsilon_n^2)-\exp(\lambda'/\varepsilon_n^2)) \geq
\exp((M(c_{m})+\delta')/\varepsilon_n^2)$ (which is true if
$\delta'$ is sufficiently small and $n$ is sufficiently large),
then from (\ref{hn}) and (\ref{hn2}) it follows that the right
hand side can be made larger than $\inf_{\overline{c} \in
[c_{m}-\delta', c_{m}+ \delta']}G(\overline{c}) - 2\delta'$. This
leads to a contradiction if $\delta'$ is small enough since $G$ is
a piece-wise continuous function which stays above the diagonal on
$[c_0 - 2\delta, c(\lambda) - \delta/2]$.

We have thus established that
$u^{\varepsilon_n}(\exp(\lambda'/\varepsilon_n^2), x_0) \leq
\beta_{i-1}$. We can then extract a further subsequence such that
 $u^{\varepsilon_n}(\exp(\lambda'/\varepsilon_n^2), x_0)$ belongs
to one of the intervals $[\beta_{j-1},\beta_j]$ with $j < i$. We
can then take $\lambda'' < \lambda'$ and repeat the argument above
to show that $u^{\varepsilon_n}(\exp(\lambda''/\varepsilon_n^2),
x_0) \leq \beta_{j-1}$. After at most $k$ such steps, we obtain
$\widetilde{\lambda} < \lambda$ such that
$u^{\varepsilon_n}(\exp(\widetilde{\lambda}/\varepsilon_n^2), x_0)
\leq \beta_{0}$ along a subsequence, and $\widetilde{\lambda}$ can
be chosen to be arbitrarily close to $\lambda$. This, however, is
a contradiction with (\ref{spll}). \qed
\\

 \noindent {\it Proof of Theorem \ref{mte}}. By
Lemma~\ref{letwo}, it is sufficient to prove that
\begin{equation} \label{center}
\lim_{\varepsilon \downarrow 0} u^\varepsilon(\exp({\lambda
/\varepsilon^{2}}), x_0) = c(\lambda).
\end{equation}

 {\it Case 1}:  $0 < \lambda
< M(c_0)$.
 Assume that
\begin{equation} \label{fs1x}
\limsup_{\varepsilon \downarrow 0} u^\varepsilon(\exp({\lambda
/\varepsilon^{2}}), x_0) > c_0.
\end{equation}
Take $c_0 < \beta_1 < \beta_2 < \limsup_{\varepsilon \downarrow 0}
u^\varepsilon(\exp({\lambda /\varepsilon^{2}}), x_0)$ such that
$M(\beta_2) > \lambda$. By Lemma~\ref{simp}, there are sequences
$\varepsilon_n \downarrow 0$ and  $v_0 \leq \mu_n < \lambda_n \leq
\lambda$ such that (\ref{ml}) holds. Thus Part (A) of
Lemma~\ref{tet} leads to a contradiction with (\ref{fs1x}). The
inequality
\[
\liminf_{\varepsilon \downarrow 0} u^\varepsilon(\exp({\lambda
/\varepsilon^{2}}), x_0) < c_0
\]
can be ruled out in the same way by referring to Part (A$'$) of
Lemma~\ref{tet}.

Case 2: $\lambda \geq M(c_0)$, $c_1 = c_0$. Assume that
\begin{equation} \label{fs1x3}
\limsup_{\varepsilon \downarrow 0} u^\varepsilon(\exp({\lambda
/\varepsilon^{2}}), x_0) > c_0.
\end{equation}
Then, since $G$ is piece-wise continuous and passes from the left
of the diagonal to the right of the diagonal at $c_1$,  we can
find $\delta > 0$ and $\beta_1$,$\beta_2$ such that
\[
c_0 < \beta_1 < \beta_2 < \limsup_{\varepsilon \downarrow 0}
u^\varepsilon(\exp({\lambda/\varepsilon^{2}}), x_0)
\]
and $G(c) < \beta_2 - \delta$ for $c \in
  [\beta_2 -\delta, \beta_2+\delta]$. By Lemma~\ref{simp}, there are sequences
$\varepsilon_n \downarrow 0$ and  $v_0 \leq \mu_n < \lambda_n \leq
\lambda$ such that (\ref{ml}) holds. Thus Part (B) of
Lemma~\ref{tet} leads to a contradiction with (\ref{fs1x3}).  The
inequality
\[
\liminf_{\varepsilon \downarrow 0} u^\varepsilon(\exp({\lambda
/\varepsilon^{2}}), x_0) < c_0
\]
can be ruled out in the same way by referring to Part (B$'$) of
Lemma~\ref{tet}.

Case 3: $\lambda \geq M(c_0)$,  $c_1
> c_0$. First assume that
\begin{equation} \label{fs1x3z}
\limsup_{\varepsilon \downarrow 0} u^\varepsilon(\exp({\lambda
/\varepsilon^{2}}), x_0) > c(\lambda).
\end{equation}
We can repeat the arguments of Case 2 to show that (\ref{fs1x3z})
implies that $c(\lambda) < c_1$. Then, since $\lambda$ is a point
of continuity of $c(\lambda)$,  we can find $\beta_1,\beta_2$ such
that
\[
c(\lambda)  < \beta_1 < \beta_2 < \limsup_{\varepsilon \downarrow
0} u^\varepsilon(\exp({\lambda /\varepsilon^{2}}), x_0)
\]
and  $M(\beta_2) > \lambda$. By Lemma~\ref{simp}, there are
sequences $\varepsilon_n \downarrow 0$ and  $v_0 \leq \mu_n <
\lambda_n \leq \lambda$ such that (\ref{ml}) holds. Thus Part (A)
of Lemma~\ref{tet} leads to a contradiction with (\ref{fs1x3z}).

Finally, from Lemma~\ref{rone} it follows that
\[
\liminf_{\varepsilon \downarrow 0} u^\varepsilon(\exp({\lambda
/\varepsilon^{2}}), x_0) \geq c(\lambda).
\]

Case 4: $\lambda \geq M(c_0)$,  $c_1 < c_0$. This is completely
similar to Case 3. \qed
\\

\noindent {\bf Remark.} If  instead of the constant $\lambda$  in
the argument of the function $u^\varepsilon$ in Theorem~\ref{mte},
we have a positive function $\lambda(\varepsilon)$ such that
$\lim_{\varepsilon \downarrow 0} \lambda(\varepsilon) = \lambda
> 0$, then
\[
\lim_{\varepsilon \downarrow 0}
u^\varepsilon(\exp({\lambda(\varepsilon) /\varepsilon^{2}}), x) =
c(\lambda).
\]
The proof of this statement requires only simple modifications to
the proof of Theorem~\ref{mte}.

\section{Exit From the Domain} \label{exitprob}
Let the differential operator $L^\varepsilon$ and the function $G$
satisfy the assumptions of Section~\ref{asc11}. Let $x \in D$ and
$\lambda > 0$. Recall that $X^{\lambda,x,\varepsilon}_s$, $s \in
[0,\exp({\lambda /\varepsilon^{2}})]$, is the process defined in
(\ref{pro}), with $\tau^\varepsilon$ being the first time when
this process reaches the boundary of $D$. We put $\tau^\varepsilon
= \infty$ on the event that the process does not reach the
boundary by the time $\exp({\lambda /\varepsilon^{2}})$. Let
$\overline{\tau}^\varepsilon = \min(\tau^\varepsilon,
\exp({\lambda /\varepsilon^{2}}))$. Thus, if $\tau^\varepsilon <
\infty$, then
$X^{\lambda,x,\varepsilon}_{\overline{\tau}^\varepsilon}$ is the
location where the process first exits the domain. Let
$\rho^\varepsilon$ be the measure on $\overline{D}$ induced by
$X^{\lambda,x,\varepsilon}_{\overline{\tau}^\varepsilon}$:
\begin{equation} \label{fox1}
\rho^\varepsilon(A) =
\mathrm{P}(X^{\lambda,x,\varepsilon}_{\overline{\tau}^\varepsilon}
\in A),~~A \in \mathcal{B}(\overline{D}).
\end{equation}
 Let
$\mu^\varepsilon$ be the restriction of $\rho^\varepsilon$ to
$\partial D$:
\begin{equation} \label{fox2}
\mu^\varepsilon(A) =
\mathrm{P}(X^{\lambda,x,\varepsilon}_{\overline{\tau}^\varepsilon}
\in A),~~A \in \mathcal{B}(\partial D).
\end{equation}
 Note that
$\mu^\varepsilon$ is not a probability measure, since
$\mathrm{P}(X^{\lambda,x,\varepsilon}_{\overline{\tau}^\varepsilon}
\in \partial D) < 1$. In this section we shall examine the
asymptotics of $\rho^\varepsilon$ and $\mu^\varepsilon$ when
$\varepsilon \downarrow 0$.

We shall distinguish several cases corresponding to different
values of~$\lambda$. First consider the case when $0 < \lambda <
M(c_0)$.
\begin{lemma} \label{fox}
If $x \in D$ and $0 < \lambda < M(c_0)$, then $\lim_{\varepsilon
\downarrow 0} \mathrm{P}(\tau^\varepsilon \leq \exp({\lambda
/\varepsilon^{2}})) =0$.
\end{lemma}
\proof From Lemma~\ref{sec2} it follows that for $\delta > 0$
there is $0 < v_0 < \lambda$ such that
\begin{equation} \label{smmm}
\lim_{\varepsilon \downarrow 0} \mathrm{P}(X^{v_0,x,\varepsilon}_s
\in D~~{\rm for}~{\rm all}~0 \leq s \leq \exp({v_0
/\varepsilon^{2}})) = 1
\end{equation}
uniformly in $x \in D^\delta$. We claim that for each $\varkappa >
0$,
\begin{equation} \label{ucno}
\lim_{\varepsilon \downarrow 0} u^\varepsilon(t,x) = c_0~~ {\rm
uniformly}~{\rm in}~ (t,x) \in [\exp({v_0 /\varepsilon^{2}}),
\exp({\lambda /\varepsilon^{2}})] \times D^\varkappa.
\end{equation}
Indeed, otherwise by Lemma~\ref{letwo} we could find sequences
$\varepsilon_n \downarrow 0$ and $\lambda_n \in (v_0, \lambda)$
such that either $\limsup_{n \rightarrow \infty}
u^{\varepsilon_n}(\exp({\lambda_n /\varepsilon_n^{2}}),x_0) > c_0$
or $\liminf_{n \rightarrow \infty}
u^{\varepsilon_n}(\exp({\lambda_n /\varepsilon_n^{2}}),x_0) <
c_0$. Suppose that the former is the case and that $c_1 = c_0$
(the argument in the cases when $c_1 > c_0$ and $c_1 < c_0$ is
similar). Then, from the conditions imposed on the function $G$ in
Section~\ref{asc11} it follows that there are $c_0 < \beta_1 <
\beta_2 < \limsup_{n \rightarrow \infty}
u^{\varepsilon_n}(\exp({\lambda_n /\varepsilon_n^{2}}),x_0)$ such
that the graph of $G$ goes below the diagonal in a neighborhood of
the interval $[\beta_1,\beta_2]$. Moreover, there are $v_0 <
\mu'_n < \lambda'_n < \lambda_n$ such that
$u^{\varepsilon_n}(\exp({\mu'_n /\varepsilon_n^{2}}),x_0) =
\beta_1$ and $u^{\varepsilon_n}(\exp({\lambda'_n
/\varepsilon_n^{2}}),x_0) = \beta_2$. This contradicts Part~(B) of
Lemma~\ref{tet}, thus establishing~(\ref{ucno}).

From (\ref{ucno}) and Part (A) of Lemma~\ref{au1} it follows that
$ \mathrm{P}( \tau^{\varepsilon} \leq \exp({\lambda
/\varepsilon^{2}})-\exp({v_0 /\varepsilon^{2}}))$ tends to zero as
$\varepsilon \downarrow 0$. From Corollary~\ref{co1} it then
follows that $\mathrm{P}(X^{\lambda,x,\varepsilon}_{\exp({\lambda
/\varepsilon^{2}})-\exp({v_0 /\varepsilon^{2}})} \in D^\delta)$
tends to one as $\varepsilon \downarrow 0$. By the Markov property
of the process and due to  (\ref{smmm}), this implies  the
statement of the lemma. \qed

Next, let us examine the case when $\lambda > \lambda_{\rm max}$.
(Recall that $\lambda_{\rm max} = \sup_{c \in [c_0,c_1]} M(c)$ if
$c_1 \geq c_0$ and $\lambda_{\rm max} = \sup_{c \in [c_1,c_0]}
M(c)$ if $c_1 \leq c_0$.)
\begin{lemma} \label{fttw}
Suppose that $x \in D$ and $\lambda
> \lambda_{\rm max}$.  If the minimum of the quasi-potential  $\min_{x \in
\partial D} {V}^{ a(\cdot, c_1)}(x_0,x)$ is achieved at a single
point $x^*(c_1)$, then $\mu^\varepsilon$ weakly converges to a
probability measure $\mu$ concentrated at $x^*(c_1)$. If the
minimum is achieved at two points $x^*_1(c_1)$ and $x^*_2(c_1)$
and $G_1(c_1) \neq G_2(c_1)$, then $\mu^\varepsilon$ weakly
converges to a probability measure $\mu$ concentrated at those two
points. In this case $\mu(x^*_1(c_1)) G_1(c_1) + \mu(x^*_2(c_1))
G_2(c_1) = c_1$.
\end{lemma}
\proof Let $\lambda_{\rm max} < \lambda' < \lambda$. Similarly to
the proof of Lemma~\ref{fox}, and using the fact that
$u^\varepsilon(\exp({\lambda /\varepsilon^{2}}),x_0)$ converges to
$c_1$ for each $\lambda > \lambda_{\rm max}$, we can show that
\begin{equation} \label{ucnov}
\lim_{\varepsilon \downarrow 0} u^\varepsilon(t,x) = c_1~~ {\rm
uniformly}~{\rm in}~ (t,x) \in [\exp({\lambda' /\varepsilon^{2}}),
\exp({\lambda /\varepsilon^{2}})] \times D^\varkappa,
\end{equation}
for each $\varkappa > 0$. Let $ \mathcal{A} = \{ x^*(c_1) \}$ if
the minimum $\min_{x \in
\partial D} {V}^{ a(\cdot, c_1)}(x_0,x)$ is achieved at a single
point, and $\mathcal{A} = \{ x^*_1(c_1), x^*_2(c_1) \}$ if the
minimum is achieved at two points. Recall that $\tau^\varepsilon$
is the first time when the process $X^{\lambda,x,\varepsilon}_s$
reaches the boundary of $D$. From (\ref{ucnov}) and
Lemma~\ref{au1} it follows that
\begin{equation} \label{jjj1}
\lim_{\varepsilon \downarrow 0} \mathrm{P}(\tau^\varepsilon <
\exp({\lambda /\varepsilon^{2}}) -  \exp({\lambda'
/\varepsilon^{2}})) = 1
\end{equation}
and for every $\delta > 0$ we have
\begin{equation} \label{jjj2}
\lim_{\varepsilon \downarrow 0} \mathrm{P}({\rm dist} (
X^{\lambda,x,\varepsilon}_{\tau^\varepsilon \wedge \exp({\lambda
/\varepsilon^{2}}) -  \exp({\lambda' /\varepsilon^{2}}) },
\mathcal{A}) \leq \delta) = 1.
\end{equation}
This immediately implies the desired result for the case of a
single minimum point.

Let $U^{\delta}(y) \subseteq \partial D$ denote the $\delta$
neighborhood of a point $y$ on the boundary. In the case when the
minimum is achieved at two points, we note that
\[
u^{\varepsilon}(\exp({\lambda /\varepsilon^{2}}),x) = \mathrm{E}
u^{\varepsilon}(\exp({\lambda
/\varepsilon^{2}})-{\tau}^{\varepsilon}\wedge (\exp({\lambda
/\varepsilon^{2}}) -  \exp({\lambda' /\varepsilon^{2}})) ,
X^{\lambda,x,\varepsilon}_{\tau^\varepsilon \wedge (\exp({\lambda
/\varepsilon^{2}}) -  \exp({\lambda' /\varepsilon^{2}}) )}),
\]
where the left hand side tends to $c_1$, while the right hand side
is equal to
\begin{equation} \label{ji1}
\mu^\varepsilon(U^{\delta}(x^*_1(c_1))) g(x^*_1(c_1)) +
\mu^\varepsilon(U^{\delta}(x^*_2(c_1))) g(x^*_2(c_1)) +
\alpha(\varepsilon),
\end{equation}
where $\lim_{\varepsilon \downarrow 0} \alpha(\varepsilon) = 0$,
as follows from (\ref{jjj1}) and (\ref{jjj2}). It also follows
from (\ref{jjj1}) and (\ref{jjj2}) that $\lim_{\varepsilon
\downarrow 0} \mu^\varepsilon(\partial D \setminus (
U^{\delta}(x^*_1(c_1)) \cup U^{\delta}(x^*_2(c_1)))) = 0$, which,
together with (\ref{ji1}), implies the desired result. \qed

Finally, we consider the case when $c_1 \neq c_0$ and the function
$c(\lambda)$ is continuous at a point $\lambda \in (M(c_0),
\lambda_{\rm max})$. The cases $c_1 < c_0$ and $c_1 > c_0$ are
completely similar to each other, so we shall only deal with the
latter one. Let us introduce the needed notations. Fix $v_0 \in
(0, c_0)$. By Theorem~\ref{mte}, $\lim_{\varepsilon \downarrow 0}
u^\varepsilon(\exp({v_0 /\varepsilon^{2}}), x_0) = c_0$  and
$\lim_{\varepsilon \downarrow 0} u^\varepsilon(\exp({\lambda
/\varepsilon^{2}}), x_0) = c(\lambda)
> c_0$. For each $c \in [c_0, c(\lambda)]$ we define
\[
\lambda^\varepsilon(c)  = \min(\inf\{\lambda' \geq v_0:
u^\varepsilon(\exp({\lambda' /\varepsilon^{2}}), x_0) \geq c\},
\lambda).
\]
Let $\alpha^\varepsilon(c)$, $c \in [c_0, c(\lambda)]$, be the
probability that the process $X^{\lambda,x,\varepsilon}_s$ reaches
the boundary of $D$ by the time $\exp({\lambda /\varepsilon^{2}})
- \exp({\lambda^\varepsilon(c) /\varepsilon^{2}})$, that is
\[
\alpha^\varepsilon(c)= \mathrm{P}(\tau^\varepsilon \leq
\exp({\lambda /\varepsilon^{2}}) - \exp({\lambda^\varepsilon(c)
/\varepsilon^{2}})).
\]
Since $\alpha^\varepsilon$ is left-continuous, it defines a
measure $\nu^\varepsilon$ on $ \mathcal{B}([c_0, c(\lambda)])$ via
$\nu^\varepsilon([c, c(\lambda)]) = \alpha^\varepsilon(c)$. It
will be important to identify the limit of $\alpha^\varepsilon$ as
$\varepsilon \downarrow 0$. We define the function  $\alpha: [c_0,
c(\lambda)] \rightarrow [0,1]$ by:
\[
\alpha(c) = 1 - \exp(\int_c^{c(\lambda)} \frac{d z}{z - G(z)}).
\]
Since $G$ is piece-wise continuous and its graph is above the
diagonal in a neighborhood of $[c_0, c(\lambda)]$, the function
$\alpha$ is a unique continuous function which satisfies the
differential equation
\begin{equation} \label{inequ}
 \alpha'(c) =   \frac{\alpha(c)-1
}{G(c) - c}
\end{equation}
in the points of continuity of $G$ and the terminal condition
$\alpha(c(\lambda)) = 0$.
%
%
Notice that $\alpha(c) \in [0,1)$ for $c \in [c_0, c(\lambda)]$.
The function $\alpha$ defines a measure $\nu$ on
$\mathcal{B}([c_0, c(\lambda)])$ via $\nu([c, c(\lambda)]) =
\alpha(c)$.
\begin{lemma}
If $c_1 > c_0$, $c(\lambda)$ is continuous at a point $\lambda \in
( M(c_0),\lambda_{\rm max})$, then
\[
\lim_{\varepsilon \downarrow 0} \alpha^\varepsilon(c) = \alpha(c)
\]
for $c \in [c_0, c(\lambda)]$.
\end{lemma}
\proof Assume first that the minimum of the quasi-potential is
achieved at a unique point $x^*(c)$ for each $c$  in a
neighborhood of $ [c_0, c(\lambda)]$, and therefore $G$ is a
continuous function there. Take $\delta, \delta'
> 0$ which will be specified later and $\delta'' > 0$ such that $x, x_0 \in D^{\delta''}$. Let
$\beta_0,\beta_1,...,\beta_k$ be such that $\beta_0 = c_0$,
$\beta_k = c(\lambda)$ and  $0 < \beta_{i} - \beta_{i-1} <
\delta$, $1 \leq i \leq k$.

Consider the processes $Y^{i, x, \varepsilon}_s =
X^{\lambda^\varepsilon(\beta_i),x,\varepsilon}_s$, $1 \leq i \leq
k-1$, and $Y^{k,x,\varepsilon}_s = X^{\lambda,x,\varepsilon}_s$.
Let $\tau^{i,\varepsilon}$ be the first time when the process
$Y^{i, x, \varepsilon}_s$ reaches the boundary of $D$. Let
$B^{i,x,\varepsilon}$, $1 \leq i \leq k-1$, be the event that
$\tau^{i,\varepsilon} \leq \exp({\lambda^\varepsilon(\beta_i)
/\varepsilon^{2}}) - \exp({\lambda^\varepsilon(\beta_{i-1})
/\varepsilon^{2}})$, and $B^{k,x,\varepsilon}$ be the event that
$\tau^{k,\varepsilon} \leq \exp({\lambda /\varepsilon^{2}}) -
\exp({\lambda^\varepsilon(\beta_{k-1}) /\varepsilon^{2}})$.

 Using
Lemma~\ref{au1}, it is not difficult to show that for each
$\delta'>0$ there is $\delta > 0$ such that
\begin{equation} \label{exe}
\lim_{\varepsilon \downarrow 0} \mathrm{P}(B^{i,x,\varepsilon}
\cap \{ {\rm dist}(Y^{i, x, \varepsilon}_{\tau^{i,\varepsilon}},
x^*(\beta_i)) \geq \delta' \}) = 0
\end{equation}
uniformly in $x \in D^{\delta''}$. Since $G$ is continuous, we can
also make sure that $\delta$ is small enough so that
\begin{equation} \label{miop}
\lim_{\varepsilon \downarrow 0} \mathrm{P}(B^{i,x,\varepsilon}
\cap \{ |g(Y^{i, x, \varepsilon}_{\tau^{i,\varepsilon}})-
G(\beta_i)| \geq \delta' \}) = 0
\end{equation}
uniformly in $x \in D^{\delta''}$. We can write $u^\varepsilon(
\exp({\lambda /\varepsilon^{2}}),x)$ in two different ways
\[
u^{\varepsilon}(\exp({\lambda /\varepsilon^{2}}),x) =
\]
\[
\mathrm{E} u^{\varepsilon}(\exp({\lambda
/\varepsilon^{2}})-{\tau}^{\varepsilon}\wedge (\exp({\lambda
/\varepsilon^{2}}) -  \exp({\lambda^\varepsilon(\beta_i)
/\varepsilon^{2}})) , X^{\lambda,x,\varepsilon}_{\tau^\varepsilon
\wedge (\exp({\lambda /\varepsilon^{2}}) -
\exp({\lambda^\varepsilon(\beta_i) /\varepsilon^{2}})) })
\]
and
\[
u^{\varepsilon}(\exp({\lambda /\varepsilon^{2}}),x) =
\]
\[
\mathrm{E} u^{\varepsilon}(\exp({\lambda
/\varepsilon^{2}})-{\tau}^{\varepsilon}\wedge (\exp({\lambda
/\varepsilon^{2}}) -  \exp({\lambda^\varepsilon(\beta_{i-1})
/\varepsilon^{2}})) , X^{\lambda,x,\varepsilon}_{\tau^\varepsilon
\wedge (\exp({\lambda /\varepsilon^{2}}) -
\exp({\lambda^\varepsilon(\beta_{i-1}) /\varepsilon^{2}})) }).
\]
Upon subtracting the right hand sides of these two equalities,
using the Markov property, Corollary~\ref{co1}, Lemma~\ref{letwo}
and (\ref{miop}), we obtain
\begin{equation} \label{nmmn}
\beta_i(1- \alpha^\varepsilon(\beta_i)) =
(\alpha^\varepsilon(\beta_{i-1}) -
\alpha^\varepsilon(\beta_i))(G(\beta_i) + h_1(i, \varepsilon))+
\beta_{i-1}(1- \alpha^\varepsilon(\beta_{i-1})) +
h_2(i,\varepsilon),
\end{equation}
where $h_1(i, \varepsilon) \leq \delta'$ and $\lim_{\varepsilon
\downarrow 0} h_2(i, \varepsilon) = 0$. This implies the desired
result once we recall that $\delta'$ and $\delta$ can be taken
arbitrarily small, since (\ref{nmmn}) shows that
$\alpha^\varepsilon$ is a type of Euler's method approximation to
the solution of (\ref{inequ}).

The condition of continuity of $G$ can be easily removed once we
recall that $G$ may have at most finitely many points of
discontinuity. \qed
\\

\noindent {\bf Remark.} Using similar arguments it is not
difficult to show that
\[
\lim_{\varepsilon \downarrow 0} \mathrm{P}( \exp({\lambda
/\varepsilon^{2}}) - \exp({\lambda^\varepsilon(c_0)
/\varepsilon^{2}}) < \tau^\varepsilon < \infty )  = 0.
\]
Moreover, using (\ref{nmmn}) and (\ref{exe}) it is possible to
show that in order to find the limit of~$\mu^\varepsilon$, one can
take $\nu$, which is the limit of $\nu^\varepsilon$, and then take
its push-forward by the function $x^*$ (since $\nu$ is an
absolutely continuous measure, it is not essential that $x^*$ may
be undefined  in a finite number of points). The push-forward of
$\nu$ will be denoted by $\mu$. Thus $\mu(A) = \nu(c \in [c_0,
c(\lambda)]: x^*(c) \in A)$, $A \in \mathcal{B}(\partial D)$.
\\

Combining this with Lemmas~\ref{fox} and \ref{fttw} and
Corollary~\ref{co1}, we can can formulate the following theorem.
\begin{theorem} \label{mte2}
Let $\rho^\varepsilon$ and $\mu^\varepsilon$ be defined by
(\ref{fox1}) and (\ref{fox2}), respectively.  If $x \in D$ and $0
< \lambda < M(c_0)$, then $\rho^\varepsilon \rightarrow
\delta_{x_0}$, where $\delta_{x_0}$ is the probability measure
concentrated at $x_0$ and $\mu^\varepsilon \rightarrow \mu$, where
$\mu$ is the trivial measure, that is $\mu(\partial D) = 0$.

If $\lambda > \lambda_{\rm max}$ and the minimum of the
quasi-potential $\min_{x \in \partial D} {V}^{ a(\cdot,
c_1)}(x_0,x)$ is achieved at a single point $x^*(c_1)$, then
$\rho^\varepsilon$ and $\mu^\varepsilon$ weakly converge to a
probability measure $\mu$ concentrated at $x^*(c_1)$. If the
minimum is achieved at two points $x^*_1(c_1)$ and $x^*_2(c_1)$
and $G_1(c_1) \neq G_2(c_1)$, then $\rho^\varepsilon$ and
$\mu^\varepsilon$ weakly converge to a probability measure $\mu$
concentrated at those two points. Moreover, in this case
$\mu(x^*_1(c_1)) G_1(c_1) + \mu(x^*_2(c_1)) G_2(c_1) = c_1$.

If $c_1 > c_0$ and the function $c(\lambda)$ is continuous at a
point $\lambda \in (M(c_0),\lambda_{\rm max})$, then take the
measure $\nu$ on $ \mathcal{B}([c_0, c(\lambda)])$ defined via
$\nu([c, c(\lambda)]) = \alpha(c)$, where $\alpha$ is the solution
of (\ref{inequ}). The measures $\mu^\varepsilon$ weakly converge
to the measure $\mu$ which is the push-forward of $\nu$ by the
function $x^*$. The measures $\rho^\varepsilon$ weakly converge to
the measure $c \delta_{x_0} + \mu$, where $c = 1 - \mu(\partial
D)$.
\end{theorem}
\noindent {\bf Remark.}  This theorem still holds if instead of
$\lambda$ in the definition of $\mu^\varepsilon$ we have a
positive function $\lambda(\varepsilon)$ such that
$\lim_{\varepsilon \downarrow 0} \lambda(\varepsilon) = \lambda>
0$.
%
\begin{corollary} If $c_1 > c_0$ and the function $c(\lambda)$ is continuous at
 $\lambda \in (M(c_0),\lambda_{\rm max})$, then for every
$\delta > 0$ and $x \in D^\delta$ we have
\[
\lim_{\varepsilon \downarrow 0} u^\varepsilon(\exp({\lambda
/\varepsilon^{2}}), x) = \int_{c_0}^{c(\lambda)} g(x^*(c)) d
\nu(c) + g(x_0) (1 - \nu([c_0, c(\lambda)]),
\]
where $\nu$ is the measure on $\mathcal{B}([c_0, c(\lambda)])$
defined via $\nu([c, c(\lambda)]) = \alpha(c)$.
\end{corollary}
\proof The corollary immediately follows from Theorems~\ref{mte}
and \ref{mte2} and the probabilistic representation of the
solution to the initial-boundary value problem. \qed

\section{Generalizations and Examples} \label{genex}
\subsection{The Case of a Nonlinear First Order Term} \label{fot}
We could allow the coefficient at the first order term to depend
on $u^\varepsilon$ in  (\ref{eq1})-(\ref{bc}):
\[
\frac{\partial u^\varepsilon(t,x)}{\partial t}  = L^\varepsilon
u^\varepsilon :=
\]
\begin{equation} \label{eq1x1}
 \frac{\varepsilon^2}{2} \sum_{i,j=1}^d
a_{ij}(x,u^\varepsilon) \frac{ \partial^2 u^\varepsilon
(t,x)}{\partial x_i \partial x_j} + (b(x, u^\varepsilon) +
\varepsilon b_1(x, u^\varepsilon)) \cdot \nabla_x u^\varepsilon
(t,x),~~x \in D,~~t
> 0,
\end{equation}
\begin{equation} \label{bcx1}
u^\varepsilon(0,x) = g(x),~x \in D,~~~~~u^\varepsilon(t,x) =
g(x),~~t \geq 0,~x
\in
\partial D.
\end{equation}
All the assumptions made in Sections~\ref{qle} and \ref{asc11}
remain in force, other than the following: instead of assuming
that $b$ is a vector valued function on $D$, we assume that $b,
b_1 \in C^2(D \times \mathbb{R})$, and there is a positive
constant $k'$ such that $(b(x,u), n(x)) < -k'$ for $x
\in
\partial D$, $u \in \mathbb{R}$, where $n(x)$ is the exterior normal to the boundary
of $D$. Moreover, we assume that for each $u$ the vector field
$b(\cdot, u)$ has a unique equilibrium point $x_0$ which does not
depend on $u$ and that all the trajectories of the dynamical
system $x'(t) = b(x(t),u)$ starting in $D$ are attracted to $x_0$.
We now assume that there is a smooth function $v$ defined on
$\overline{D}$, such that $v(x_0) = 0$, $v(x) > 0$ for $x \neq
x_0$, and $(b(x,u), \nabla v(x)) \leq -c|x-x_0|^2$ for some
positive constant $c$, all $u$ and all $x$.

The definition of the function $M(c)$ from Section~\ref{asc11}
needs to be modified to allow for the dependence of the drift term
on a parameter. Namely, now
\[
 M(c) = \min_{x \in
\partial D} {V}^{ a(\cdot, c), b(\cdot, c)}(x_0,x),
\]
where ${V}^{ a(\cdot, c), b(\cdot, c)}$ is the quasi-potential for
the process whose generator is equal to
\[
 \frac{\varepsilon^2}{2}
\sum_{i,j=1}^d a_{ij}(x,c) \frac{ \partial^2 u^\varepsilon
(t,x)}{\partial x_i \partial x_j} + b(x,c) \cdot\nabla_x
u^\varepsilon (t,x).
\]
With this definition of $M(c)$, Theorems~\ref{mte} and \ref{mte2}
remain valid, and the proofs do not require serious modifications.

\subsection{Metastable Distributions in the Case of Two Equilibrium Points}
In this section we again consider the solutions $u^\varepsilon$ to
(\ref{eq1})-(\ref{bc}). Let all the assumptions about the domain
$D$ and the operator $L^\varepsilon$ made in Sections~\ref{qle}
and \ref{fot} remain in force, except the following: instead of
assuming the existence of a singe equilibrium point, we assume
that there are two asymptotically stable equilibrium points $x_1,
x_2 \in D$ such that for almost every $x \in D$, with respect to
the Lebesgue measure, the trajectory  of (\ref{dsyst}) starting at
$x$ is attracted to either  $x_1$ or $x_2$. Let $D_1 \subset D$ be
the set of points in $D$ which are attracted to $x_1$ and $D_2
\subset D$ the set of points attracted to $x_2$.

As before, we need to study the quasi-potential in order to
determine the asymptotic behavior of $u^\varepsilon$. While in the
case of a single equilibrium, the function $u^\varepsilon$ was
nearly constant in $D^\delta$ at times of order
$\exp(\lambda/\varepsilon^2)$ (Lemma~\ref{letwo}), now
$u^\varepsilon(\exp(\lambda/\varepsilon^2), x)$ will be close to
$u^\varepsilon(\exp(\lambda/\varepsilon^2), x_1)$ for $x \in
D_1^\delta$ and close to
$u^\varepsilon(\exp(\lambda/\varepsilon^2), x_2)$ for $x \in
D_2^\delta$. This explains why instead of freezing the second
variable in the coefficients $a_{ij}$ in the right hand side of
(\ref{eq1}), the way it was done is Section~\ref{qle}, now we put
the variable equal to $c_1$ in $D_1$ and $c_2$ in $D_2$. More
precisely, for $c_1,c_2 \in [g_{\rm min}, g_{\rm max}]$, let
\[
f_{c_1, c_2}(x) = c_1 \chi_{D_1}(x) + c_2 \chi_{D_2}(x) + (c_1 +
c_2) \chi_{D \setminus (D_1 \cup D_2)} (x)/2,~~x \in D,
\]
where $\chi_U$ is the indicator function of a set $U \subseteq
\mathbb{R}^d$. For a measurable positive-definite matrix-valued
function $\alpha$ on $\overline{D}$, we define
\[
V^\alpha(x,y) = \inf_{T, \varphi} \{ S^\alpha_{0,T}(\varphi):
\varphi \in C([0,T], \overline{D}),  \varphi(0) = x, \varphi(T) =
y \},~~x,y \in \overline{D},
\]
where the normalized action functional $S$ was defined in
Section~\ref{af}. Instead of function $M$ used in
Section~\ref{qle}, we now have functions $M_{x_1,x_2}$,
$M_{x_2,x_1}$, $M_{x_1, \partial D}$ and $M_{x_2, \partial D}$.
These are defined by
\begin{equation} \label{eqTxx1}
M_{x_1,x_2}(c_1) =  {V}^{ a(\cdot, f_{c_1, c_2}(\cdot))}(x_1,x_2),
\end{equation}
\begin{equation} \label{eqTxx2}
M_{x_2,x_1}(c_2) =  {V}^{ a(\cdot, f_{c_1, c_2}(\cdot))}(x_2,x_1),
\end{equation}
\[
M_{x_1, \partial D}(c_1,c_2) = \min_{x \in
\partial D} {V}^{ a(\cdot, f_{c_1, c_2}(\cdot))}(x_1,x),
\]
\[
M_{x_2, \partial D}(c_1,c_2) = \min_{x \in
\partial D} {V}^{ a(\cdot, f_{c_1, c_2}(\cdot))}(x_2,x).
\]
It is not difficult to check that the right hand side of
(\ref{eqTxx1}) does not depend on $c_2$ and the right hand side of
(\ref{eqTxx2}) does not depend on $c_1$. For the process governed
by equation~(\ref{eq1}), with $a_{ij}(\cdot, f_{c_1, c_2}(\cdot))$
instead of the nonlinear coefficients $a_{ij}(x,u^\varepsilon)$,
the transition from $x_1$ to a small neighborhood of $x_2$ occurs
in time of order $\exp(M_{x_1,x_2}(c_1)/\varepsilon^2)$ (provided
that the process does not exit the domain $D$ earlier). Similarly,
the transition from $x_2$ to a neighborhood of $x_1$ occurs in
time of order $\exp(M_{x_2,x_1}(c_2)/\varepsilon^2)$, while the
transition from $x_1$ and $x_2$ to the boundary occurs in time of
order $\exp(M_{x_1,\partial D}(c_1,c_2)/\varepsilon^2)$ and
$\exp(M_{x_2,\partial D}(c_1,c_2)/\varepsilon^2)$, respectively.

In this section we would like to study the equation at a time
scale which is sufficiently large for the process to make
excursions between the neighborhoods of $x_1$ and $x_2$ and back,
yet not too large so that the process starting at $x_1$ or $x_2$
does not exit the domain. Therefore we assume that
\[
\max( \sup_{c_1 \in [g_{\rm min}, g_{\rm max}]} M_{x_1,x_2}(c_1),
\sup_{c_2 \in [g_{\rm min}, g_{\rm max}]} M_{x_2,x_1}(c_2)) <
M^\partial,
\]
where
\[
M^\partial = \min( \inf_{c_1, c_2 \in [g_{\rm min}, g_{\rm max}]}
M_{x_1,
\partial D}(c_1,c_2), \inf_{c_1, c_2 \in [g_{\rm min}, g_{\rm max}]} M_{x_2,
\partial D}(c_1,c_2)).
\]
For example, if $a$ and $b$ are defined in the entire space $
\mathbb{R}^d$, $a$ is bounded and $b$ satisfies $(b(x), x) \leq A
- B(x,x)$ for some positive constants $A$ and $B$, then this
condition will be satisfied for any domain $D$ which contains a
sufficiently large ball centered at the origin.

Let $c_1 = g(x_1)$ and $c_2 = g(x_2)$. Without loss of generality
we may assume that $c_1 \leq c_2$.  In order to formulate the
results on the asymptotics of $u^\varepsilon$, we need to
introduce functions $c^1(\lambda)$ and $c^2(\lambda)$ which are
similar to the function $c(\lambda)$ from in Section~\ref{asc11}.

Let $0< \lambda < M^\partial$, and define $c^1(\lambda)$ as
follows:

For $0 < \lambda < M_{x_1,x_2}(c_1)$, let $c^1(\lambda) = c_1$.

For $\lambda \geq M_{x_1,x_2}(c_1)$, let $ c^1(\lambda) = \min
\{c_2, \min\{c: c \in [c_1, c_2], M_{x_1,x_2}(c) = \lambda \} \}$.

\noindent Similarly, we define $c^2(\lambda)$ as follows:

For $0 < \lambda < M_{x_2,x_1}(c_2)$, let $c^2(\lambda) = c_2$.

For $\lambda \geq M_{x_2,x_1}(c_2)$, let $ c^2(\lambda) =
\max\{c_1, \max\{c: c \in [c_1, c_2], M_{x_2,x_1}(c) = \lambda \}
\}$.

 Let $\lambda^* = \inf\{\lambda:
c^1(\lambda) \geq c^2(\lambda)\}$. Assume that at least one of the
functions $c^1$ and $c^2$ is continuous at $\lambda^*$. Let $c^* =
c^1(\lambda^*)$ if $c^1$ is continuous at $\lambda^*$ and $c^* =
c^2(\lambda^*)$ otherwise. Let $\overline{c}^1(\lambda) =
\min(c^1(\lambda), c^*)$ and $\overline{c}^2(\lambda) =
\max(c^2(\lambda), c^*)$.

 We can now formulate the following
analog of Theorem~\ref{mte}

\begin{theorem} \label{mtex1}
Let the above assumptions  be satisfied. Suppose that the function
$\overline{c}^1(\lambda)$ is continuous at a point $\lambda \in
(0, M^\partial)$. Then for every $\delta > 0$ the following limit
\[
\lim_{\varepsilon \downarrow 0} u^\varepsilon(\exp({\lambda
/\varepsilon^{2}}), x) = \overline{c}^1(\lambda)
\]
is uniform in $x \in D_1^\delta$. Suppose that the function
$\overline{c}^2(\lambda)$ is continuous at a point $\lambda \in
(0, M^\partial)$. Then for every $\delta > 0$ the following limit
\[
\lim_{\varepsilon \downarrow 0} u^\varepsilon(\exp({\lambda
/\varepsilon^{2}}), x) = \overline{c}^2(\lambda)
\]
is uniform in $x \in D_2^\delta$.
\end{theorem}

\noindent {\bf Remark.} If $\lambda > \lambda^*$, then $
\overline{c}^1(\lambda) = \overline{c}^2(\lambda) =  c^*$.  It is
not difficult to see that the limit
\[
\lim_{\varepsilon \downarrow 0} u^\varepsilon(\exp({\lambda
/\varepsilon^{2}}), x) = c^*
\]
is uniform in $(x, \lambda) \in D^\delta \times
[\overline{\lambda},\infty)$ for each $\overline{\lambda} >
\lambda^*$. Therefore, for each $\delta > 0$ and
$\overline{\lambda}
> \lambda^*$ there is $\varepsilon_0
> 0$ such that
\[
|u^\varepsilon(t, x) - c^*| \leq \delta
\]
whenever $\varepsilon \in (0,\varepsilon_0)$, $x \in D^\delta$ and
$t \geq \exp({\overline{\lambda} /\varepsilon^{2}})$.

Recall that $X^{\lambda,x,\varepsilon}_s$, $s \in [0,\exp({\lambda
/\varepsilon^{2}})]$, is the process defined in (\ref{pro}),
$\tau^\varepsilon$ is the first time when this process reaches the
boundary of $D$ and $\overline{\tau}^\varepsilon =
\min(\tau^\varepsilon, \exp({\lambda /\varepsilon^{2}}))$. Since
we assume that $\lambda < M^\partial$, the probability that $
\tau^\varepsilon < \exp({\lambda /\varepsilon^{2}})$ now tends to
zero as $\varepsilon \downarrow 0$. The distribution of the random
variable $X^{\lambda,x,\varepsilon}_{\overline{\tau}^\varepsilon}$
will be concentrated near the points $x_1$ and $x_2$.
\begin{theorem} \label{tenn}
Suppose that $c_1 \neq c_2$. If the function
$\overline{c}^1(\lambda)$ is continuous at a point $\lambda \in
(0, M^\partial)$ and $x \in D_1$, then the distribution of the
random variable
$X^{\lambda,x,\varepsilon}_{\overline{\tau}^\varepsilon}$
converges to the measure $\mu_1^\lambda = a_1 \delta_{x_1} + a_2
\delta_{x_2}$, where the coefficients $a_1$ and $a_2$ can be found
from the equations $\overline{c}^1(\lambda) = a_1 c_1 + a_2 c_2$,
$a_1 + a_2 = 1$.

If the function $\overline{c}^2(\lambda)$ is continuous at a point
$\lambda \in (0, M^\partial)$ and $x \in D_2$, then the
distribution of the random variable
$X^{\lambda,x,\varepsilon}_{\overline{\tau}^\varepsilon}$
converges to the measure $\mu_2^\lambda = a_1 \delta_{x_1} + a_2
\delta_{x_2}$, where the coefficients $a_1$ and $a_2$ can be found
from the equations $\overline{c}^2(\lambda) = a_1 c_1 + a_2 c_2$,
$a_1 + a_2 = 1$.

If $\lambda \in (\lambda^*, M^\partial)$ and $x \in D$, then the
distribution of the random variable
$X^{\lambda,x,\varepsilon}_{\overline{\tau}^\varepsilon}$
converges to the measure $\mu^* = a_1 \delta_{x_1} + a_2
\delta_{x_2}$, where the coefficients $a_1$ and $a_2$ can be found
from the equations $c^* = a_1 c_1 + a_2 c_2$, $a_1 + a_2 = 1$.
\end{theorem}

The proofs of Theorems~\ref{mtex1} and \ref{tenn} rely on the same
techniques as those used in the proofs of Theorems~\ref{mte} and
\ref{mte2}, and therefore will not be presented here. If $a$ is
bounded and $b$ satisfies $(b(x), x) \leq A - B(x,x)$ for some
positive constants $A$ and $B$, then similar results can be
formulated for the Cauchy problem in $ \mathbb{R}^d$ and the
corresponding nonlinear perturbations of the dynamical system in $
\mathbb{R}^d$. In this case we do not need the condition $\lambda
< M^\partial$, but can instead consider all $\lambda \in (0,
\infty)$.

\subsection{Examples} In this section we give two examples when we
can easily calculate the function $M(c)$ defined in
Section~\ref{asc11}.

In the first example, we assume that the domain $D$ is one
dimensional: $D = (A, B)$. We assume that $a(x,u) \in
C^2([A,B]\times \mathbb{R})$ is positive, $b(x,u) \in
C^2([A,B]\times \mathbb{R})$, $\partial b(x,u)/\partial x < k < 0$
and there is a point $x_0 \in (A,B)$ such that $b(x_0,u) = 0$ for
all $u$. In this case the operator
\[
L^\varepsilon u = \frac{\varepsilon^2}{2}a(x,u)u'' + b(x,u)u'
\]
satisfies the assumptions of Section~\ref{fot}. The
quasi-potential, which will now be denoted by $V^c$,  is given by
\[
{V}^c(x_0,x) = -2 \int_{x_0}^x \frac{b(y,c)}{a(y,c)} d y,
\]
as is easily seen from the definition of the action functional
(see Section~\ref{af}). Therefore,
\[
M(c) = \min(V^c(x_0,A), V^c(x_0, B)) =  \min(-2 \int_{x_0}^A
\frac{b(y,c)}{a(y,c)} d y, -2 \int_{x_0}^B \frac{b(y,c)}{a(y,c)} d
y).
\]
The function $G(c)$ may take at most two values: $g(A)$ and
$g(B)$. In particular, the value $g(A)$ is taken on the set
$\{c:V^c(x_0,A) < V^c(x_0, B) \}$, while the value $g(B)$ is taken
on the set  $\{c: V^c(x_0,A) > V^c(x_0, B)\}$.

Picture 2 shows an example of the graphs of functions $M(c)$ and
$G(c)$ in the case when $g(A) = \inf_{x \in [A, B]} g(x) < \sup_{x
\in [A, B]} g(x) = g(B)$. From Theorem~\ref{mte2} and the
discussion in Section~\ref{fot} it follows that for $\lambda >
\lambda_{\rm max}$, the distribution of the random variable
$X^{\lambda,x,\varepsilon}_{\overline{\tau}^\varepsilon}$
converges to the probability measure $\mu$ concentrated at the end
points of the segment. This measure can be found from the
relations
\[
\mu(A) g(A) + \mu(B) g(B) = c_1,~~~\mu(A) + \mu(B) = 1.
\]

\begin{figure}[htbp]
 \label{pic2}
  \begin{center}
    \begin{psfrags}
     \includegraphics[height=4.0in, width= 4.5in,angle=0]{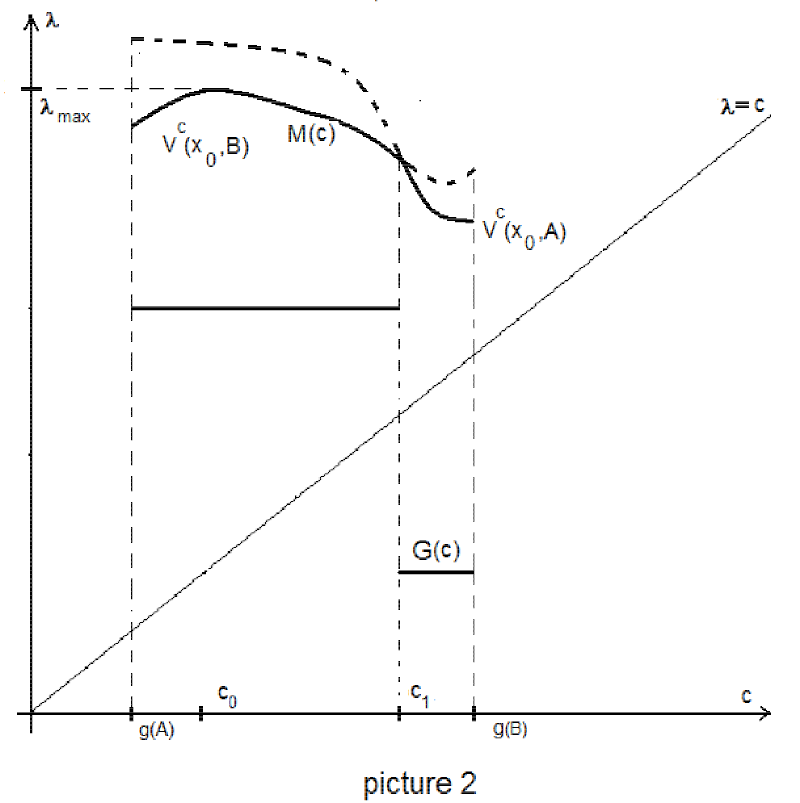}
    \end{psfrags}
  \end{center}
\end{figure}

 In the second example we assume that $D \subset \mathbb{R}^d$
contains the origin $x_0 = 0$. Let the operator $L^\varepsilon$ be
as follows
\[
L^\varepsilon u = \frac{\varepsilon^2}{2} \sum_{i,j=1}^d a_{ij}(u)
\frac{ \partial^2 u (t,x)}{\partial x_i \partial x_j} + (A(u) x)
\cdot \nabla u,
\]
where $a$ is a positive-definite matrix which depends smoothly on
$u$, and $A$ is a matrix with negative eigenvalues which depends
smoothly on $u$. As has been demonstrated in~\cite{CF}, the
quasi-potential, which we shall denote by $V^c$, is given by the
quadratic form
\[
V^c(x_0,x) = \frac{1}{2}(B^{-1}(c)x,x),
\]
where the matrix $B$ is given by
\[
B(c) = \int_0^\infty \exp(A(c)t) a(c) \exp(A^*(c) t) d t.
\]
Therefore, $M(c) = \min_{x \in \partial D} (B^{-1}(c)x,x)/2$.
\\
\\
\noindent {\bf \large Acknowledgements}: While working on this
article, M. Freidlin was supported by NSF grant DMS-0803287 and L.
Koralov was supported by NSF grant DMS-0706974.


\begin{thebibliography}{999999}

\bibitem{Ar} Aronson D.G., {\it Non-negative solutions of linear
parabolic equations}, Ann. Sci. Norm. Sup. Pisa 22 (1968),
607-694.

\bibitem{AF} Athreya A., Freidlin M.I., {\it Metastability and
Stochastic Resonance in Nearly-Hamiltonian Systems}, Stochastics
and Dynamics, 8, 1, pp 1-21, 2008.

\bibitem{CF} Chen Z., Freidlin M.I., {\it Smoluchowski-Kramers
Approximation and Exit Problems}, Stochastics and Dynamics, 5, pp
569-585, 2005.

\bibitem{F} Freidlin M. I., {\it Metastability and Stochastic
Resonance of Multiscale Systems}, Contemporary Mathematics, 2008.

\bibitem{F77} Freidlin M. I., {\it Sublimiting Distributions and
Stabilization of Solutions of Parabolic Equations with a Small
Parameter}, Soviet Math. Dokl., 235, 5, pp 1042-1045, 1977.

\bibitem{FPhD} Freidlin M. I., {\it Quasi-deterministic
Approximation, Metastability and Stochastic Resonance}, Physica D,
137, pp 333-352, 2000.

\bibitem{FK2} Freidlin M.I., Koralov L., {\it Metastability
for Nonlinear Random Perturbations of Dynamical Systems},
Stochastic Processes and Applications 120 (2010), no. 7,
1194-1214.

\bibitem{FW} Freidlin M. I., Wentzell A. D., {\it Random
Perturbations of Dynamical Systems}, Springer 1998.

\bibitem{K} Krylov  N.V., {\it  Nonlinear Elliptic and Parabolic Equations of the Second Order
 (Mathematics and its Applications)}, Springer 1987.

\bibitem{OV} Oliveiri E., Vares M.E. {\it Large Deviations and
Metastability}, Cambridge University Press, 2005.

\end{thebibliography}
\end{document}